\documentclass[a4paper,11pt, reqno]{amsart}

\usepackage{amsmath}
\usepackage{amssymb}
\usepackage{mathrsfs}
\usepackage{ifthen}
\usepackage{graphicx}
\usepackage[T1]{fontenc} 
\usepackage{lipsum}
\usepackage[a4paper,margin=0.75in]{geometry}

\usepackage{cite}

\def\switchlinenumbers{\@ifstar
	{\let\makeLineNumberOdd\makeLineNumberRight
		\let\makeLineNumberEven\makeLineNumberLeft}%
	{\let\makeLineNumberOdd\makeLineNumberLeft
		\let\makeLineNumberEven\makeLineNumberRight}%
}

\def\setmakelinenumbers#1{\@ifstar
	{\let\makeLineNumberRunning#1%
		\let\makeLineNumberOdd#1%
		\let\makeLineNumberEven#1}%
	{\ifx\c@linenumber\c@runninglinenumber
		\let\makeLineNumberRunning#1%
		\else
		\let\makeLineNumberOdd#1%
		\let\makeLineNumberEven#1%
		\fi}%
}

\nonstopmode \numberwithin{equation}{section}
\newtheorem*{theorem*}{Theorem}

\newtheorem{thm}{Theorem}[section]
\newtheorem{cor}{Corollary}[section]
\newtheorem{lem}{Lemma}[section]

\theoremstyle{definition}
\newtheorem{defn}{Definition}[section]
\newtheorem{example}{Example}[section]
\newtheorem{qsn}{Question} [section]
\newtheorem{prob}[equation]{Problem}
\newtheorem{rem}{Remark}[section]

\newtheorem{innercustomthm}{Theorem}
\newenvironment{customthm}[1]
{\renewcommand\theinnercustomthm{#1}\innercustomthm}
{\endinnercustomthm}

\newcounter{minutes}
\divide\time by 60
\newcounter{hours}
\multiply\time by 60
\addtocounter{minutes}{-\time}

\newcounter {own}
\def\theown {\thesection       .\arabic{own}}

\newenvironment{pf}[1][]{%
	\vskip 3mm
	\noindent
	\ifthenelse{\equal{#1}{}}%
	{{\slshape {\bf Proof}. }}%
	{{\slshape #1.} }%
}%
{\qed\bigskip}

\newcounter{alphabet}

\def\be{\begin{equation}}
	\def\ee{\end{equation}}

\newcommand{\bee}{\begin{enumerate}}
	\newcommand{\eee}{\end{enumerate}}

\newcommand{\blem}{\begin{lem}}
	\newcommand{\elem}{\end{lem}}
\newcommand{\bthm}{\begin{thm}}
	\newcommand{\ethm}{\end{thm}}
\newcommand{\bcor}{\begin{cor}}
	\newcommand{\ecor}{\end{cor}}
\newcommand{\beg}{\begin{examp}}
	\newcommand{\eeg}{\end{examp}}
\newcommand{\begs}{\begin{examples}}
	\newcommand{\eegs}{\end{examples}}
\newcommand{\bdefe}{\begin{defin}}
	\newcommand{\edefe}{\end{defin}}
\newcommand{\bprob}{\begin{prob}}
	\newcommand{\eprob}{\end{prob}}
\newcommand{\bei}{\begin{itemize}}
	\newcommand{\eei}{\end{itemize}}
\newcommand{\real}{{\operatorname{Re}\,}}

\newcommand{\norm}[1]{\left\lVert#1\right\rVert}

\newcommand{\comment}[1]{}

\subjclass[{AMS} Subject Classification:]{Primary 32A05, 31C10, 46B07;  Secondary 32Q02, 46E40}

\begin{document}
	
	\title[]{Arithmetic Bohr radius and Local Banach space theory
	}

	\author{Himadri Halder}
	\address{Himadri Halder,
		Department of Mathematics,
		Indian Institute of Science, Bangalore-560012, India}
	\email{himadrihalder119@gmail.com, himadrih@iisc.ac.in}
	

	
	
	\begin{abstract}
		This article introduces the notion of arithmetic Bohr radius for operator valued pluriharmonic functions on complete Reinhardt domains in $\mathbb{C}^n$. Using tools from local Banach space theory, we determine its asymptotic behavior in both finite and infinite dimensions. Asymptotic estimates for this constant are derived for both convex and non-convex complete Reinhardt domains. The framework developed in this article extends the classical Minkowski-space setting to a much broader class of sequence spaces, such as mixed Minkowski, Lorentz, and Orlicz spaces. Our results also apply to a wide class of Banach sequence spaces, including symmetric and convex Banach spaces. This generality allows for a unified and systematic investigation of Bohr’s theorem for both holomorphic and pluriharmonic functions. As an application of our results, we obtain several consequences extending known results in the scalar valued setting and in the existing literature.
	\end{abstract}

	\maketitle
	\pagestyle{myheadings}
	\markboth{Himadri Halder}{Bohr and arithmetic Bohr radii via local Banach space theory}
	
	\section{Introduction and the main results}\label{section-1}
	Over the past three decades, Bohr's theorem and its applications have become a central topic in the function theory of several complex variables, exhibiting deep connections with fundamental aspects of both the local and global theory of Banach spaces. These rich connections have stimulated extensive research on Bohr's theorem in various fields of mathematics. For instance, it has been investigated in the contexts of Banach algebras and uniform algebras (see \cite{paulsen-2002,paulsen-2004}), complex manifolds (see \cite{aizn-2000b,aizenberg-2001a}), ordinary and vector-valued Dirichlet series (see \cite{bala-2006,defant-2008}), elliptic equations (see \cite{aizenberg-2001b}), Faber–Green condensers (see \cite{lassere-2017}), free holomorphic functions (see \cite{popescu-2019}), vector-valued holomorphic functions (see \cite{defant-2011}), local Banach space theory (see \cite{defant-2003,defant-2025-TAMS}), domains of monomial convergence (see \cite{defant-2009}), harmonic and pluriharmonic mappings (see \cite{hamada-JFA-2022}), slice regular functions (see \cite{Rocchetta-Slice-Bohr-Math-Nachr,Xu-quaternion-AMPA,Xu-quaternion-PRSESA}), Hardy spaces (see \cite{bene-2004}), as well as various multidimensional settings (see \cite{aizn-2000a,boas-1997,boas-2000,defant-2004,defant-2006}). 
	\par 
	The investigation of Bohr's theorem for slice regular functions has recently attracted considerable attention in hypercomplex analysis. The pioneering work of Rocchetta, Gentili, and Sarfatti \cite{Rocchetta-Slice-Bohr-Math-Nachr} in $2012$ established the Bohr's theorem for slice regular functions, providing a fundamental extension of the classical Bohr phenomenon to the setting of noncommutative function theory. Their results opened new directions in the study of analytic properties of quaternionic and related function spaces. Subsequently, Xu \cite{Xu-quaternion-AMPA} in $2021$  made significant contributions by introducing the notion of the generalized Bohr radius for slice regular functions over quaternions, thereby refining and extending earlier results. In the same year, Xu \cite{Xu-quaternion-PRSESA} further extended the Bohr's theorem to slice regular functions over octonions, addressing additional challenges arising from nonassociativity. These developments have considerably enriched the theory and stimulated further research on Bohr's theorem in hypercomplex frameworks.
	\par In multidimensional settings, the Bohr radius problem becomes significantly more complicated. Many questions concerning the precise behavior of multidimensional Bohr radii remain open, largely due to the lack of certain essential analytical and geometric tools. Several of these problems are closely connected to fundamental topics in functional analysis, such as the geometry of Banach spaces \cite{Blasco-Collect-2017,defant-2012}, unconditional basis constants of spaces of polynomials \cite{aizenberg-2001a}, as well as Gordon–Lewis constants, projection constants, and the Banach–Mazur distance \cite{defant-2003}.
	 An especially intriguing direction arises when one considers this problem for vector valued functions on complete Reinhardt domains.
	\par	
	Our goal of this article is to build a local Banach space theoretic framework for the study of Bohr's theorem for vector valued holomorphic and pluriharmonic functions via the arithmetic Bohr radius. One of the core challenges in this study is to analyze the asymptotic behavior of both the Bohr radius and the arithmetic Bohr radius using tools and invariants from local Banach space theory. To proceed in more detail, we first recall the notion of the Bohr radius for vector valued holomorphic functions.
	Let $\Omega \subset \mathbb{C}^n$ be a given complete Reinhardt domain and $n \in \mathbb{N}$. Let $U:X\rightarrow Y$ be a non-null bounded liner operator between two complex Banach spaces $X$ and $Y$, with $\norm{U} \leq \lambda$. For $1 \leq p < \infty$, the $\lambda$-powered Bohr radius of $U$, denoted by $K_{\lambda}(\Omega, p,U)$, is defined to be the supremum of all $r\geq 0$ such that for all holomorphic functions $f(z)=\sum_{\alpha}a_{\alpha}z^{\alpha}:\Omega \rightarrow X$ we have 
	\begin{equation} \label{e-1.2}
		\sup_{z \in r\Omega}\, \sum_{\alpha} \norm{U(a_{\alpha})z^\alpha}^p_{Y} \leq \lambda^p\,\norm{f}^p_{\Omega,X}, 
	\end{equation}
	where $\norm{f}_{\Omega,X}:=\sup_{z \in \Omega}\norm{f(z)}_{X}$. 
	\vspace{2mm}
	
	Let us fix some notations. We write $K(\Omega, p,U):=K_{1}(\Omega, p,U)$. When $U=I:X\rightarrow X$ is the identity operator, we set $K_{\lambda}(\Omega, p,X):=K_{\lambda}(\Omega, p,I)$ and $K(\Omega, p,X):=K_{1}(\Omega, p,X)$. In the scalar valued case, we denote $K_{\lambda}(\Omega, p):=K_{\lambda}(\Omega, p,\mathbb{C})$ and $K(\Omega,p):=K_1(\Omega,p)$.
	\vspace{1mm}
	  
	The above definition is motivated by Defant {\it et al.} \cite{defant-2012}, extending their result for the case $p=1$ and $\Omega=\mathbb{D}^n$ to a more general setting, thus making it possible to the study of Bohr phenomenon for broader classes of domains and functions. 
	With the above notations, a classical theorem of Harald Bohr concerning power series states that $K(\mathbb{D},1)=1/3$ (see \cite{Bohr-1914}). 
	One of the most striking features of Bohr's theorem is that our understanding of the behavior of the constant $K(\Omega,1)$ remains rather limited. A central challenge is to determine the precise value of $K(\Omega,1)$, which is still unknown for $n>1$, even in the classical setting of the polydisc $\mathbb{D}^n$, where $\mathbb{D}^n:=\{z=(z_{1},\ldots,z_{n}) \in \mathbb{C}^n: |z_{j}|<1\, \mbox{for all}\, 1 \leq j \leq n\}$. To make progress in this direction, much of the existing research has focused on establishing bounds for $K(\Omega,1)$ and investigating its asymptotic behavior on general domains, and in particular well known settings such as the unit ball of Minkowski space and finite-dimensional complex Banach spaces. 
	\par
	Denote $B_{\ell^n _q}:=\left\{z=(z_{1},\ldots,z_{n}) \in \mathbb{C}^n:\norm{z}_{q}:=\left(\sum_{i=1}^{n}|z_{i}|^q\right)^{1/q}<1\right\}$, $1\leq q <\infty$ and $B_{\ell^n _\infty}:=\mathbb{D}^n$. For every $1\leq q \leq \infty$ and all $n\in \mathbb{N}$, together the results of Aizenberg, Boas, Dineen, Khavinson, Timoney, Defant, and Frerick from \cite{aizn-2000a,boas-1997,boas-2000,defant-2011-lpn} show that  there exist constants $C,D>0$ such that 
	\begin{equation*}
		\frac{1}{C}\, \left(\frac{\log\,n}{n}\right)^{1-\frac{1}{\min\{q,2\}}} \leq K(B_{\ell^n _q},1) \leq D \,\left(\frac{\log\,n}{n}\right)^{1-\frac{1}{\min\{q,2\}}}.
	\end{equation*}
Boas and Khavinson \cite{boas-1997} showed that the constant $K(\mathbb{D}^n,1)$ does not exceed $ 1/3$, and moreover that the sequence $\{K(\mathbb{D}^n,1)\}_{n=1}^{\infty}$ is decreasing with limit zero as $n \rightarrow \infty$. In the bidisc case, their work yields the bounds $0.23570226\leq K(\mathbb{D}^2,1) \leq 1/3=0.33333$. Considerable progress has been made only recently. Knese \cite{knese-2025} established in $2025$ that $K(\mathbb{D}^2,1)>0.3006$, while a subsequent result in $2026$ of Baran, Pikul, Woerdeman, and Wojtylak \cite{baran-2026} showed that $K(\mathbb{D}^2,1) < 0.3177$. These advances in particular ensure that $K(\mathbb{D}^n,1) <1/3$ for $n \geq 2$.
	\vspace{1mm}
	
    It is worth mentioning that the constant $K_{\lambda}(\Omega, p,U)$ is not necessarily nonzero for arbitrary Banach spaces $X,Y$, and the operator $U$. We refer to the occurrence of strictly positive constant $K_{\lambda}(\Omega, p,U)$ as the Bohr phenomenon. More precisely, we make the following definition. For given Banach spaces $X,Y$ and the operator $U$ as above, a class of holomorphic functions $f$, defined in a complete Reinhardt domain $\Omega$, is said to exhibit the Bohr phenomenon on $\Omega$ if there exists a universal constant $r=r_{0}\in (0,1]$, called the Bohr radius of $\Omega$ with respect to this class, such that inequality \eqref{e-1.2} is satisfied by every function in the class. Since the Bohr phenomenon does not occur for all classes of holomorphic functions (see \cite{aizn-2000b,Himadri-Vasu-PEMS,bene-2004,Blasco-Collect-2017}), an important problem is to identify and characterize those classes for which the Bohr phenomenon holds. Blasco \cite{Blasco-Collect-2017} initiated the systematic study of classes exhibiting the Bohr phenomenon by characterizing them in the classical one-dimensional setting $\Omega=\mathbb{D}$ with $U$ equal to the identity operator. Many thanks to Defant, Maestre, and Schwarting \cite{defant-2012}, who were the first to address this problem in multidimensional case $\Omega=\mathbb{D}^n$, $n>1$, with $p=1$ and for an arbitrary operator $U$. More recently, the problem for general bounded complete Reinhardt domains has been considered by the present author in \cite{Himadri-local-Banach-1}, where a framework based on local Banach space theory was developed.
    \vspace{2mm}
    
		The constant $K_{\lambda}(\Omega, 1)$ was first introduced by Bombieri \cite{bombieri-1962}, who determined its exact value for $\lambda \in [1,\sqrt{2}]$. Its precise asymptotic behavior as $\lambda\rightarrow \infty$ was later investigated by Bombieri and Bourgain \cite{bombieri-2004}. Building upon these foundational results,  Defant {\it et al.} \cite{defant-2012} have introduced the more general constant $K_\lambda(\mathbb{D}^n,1,U)$ and obtained asymptotic estimates in both finite and infinite dimensional Banach spaces $X$. On the other hand, the constants $K(\mathbb{D},p)$ and $K(\mathbb{D}^n,p)$ were first studied by Djakov and Ramanujan \cite{Djakov & Ramanujan & J. Anal & 2000}, and subsequently developed by B\'{e}n\'{e}teau, Dahlner, Khavinson, Das \cite{bene-2004}. From the viewpoint of Banach space valued functions, Blasco \cite{Blasco-Collect-2017} has shown that $K(\mathbb{D},p,X)=0$ for $p \in [1,2)$, whereas  $K(\mathbb{D},p,X)>0$ precisely when $X$ is $p$-uniformly $\mathbb{C}$-convex for $2 \geq p < \infty$. 
	These results show that $K(\mathbb{D},p,X)$ is not necessarily nonzero for all Banach spaces $X$. This observation inspired a refinement of the definition to guarantee non-vanishing behavior for arbitrary $X$, leading to the introduction of the parameter $\lambda$ in \cite[Definition 1.1]{defant-2012}, and hence in \eqref{e-1.2}.
	\vspace{2mm}
	
	One of the main aims of this note is to study the Bohr radius through the lens of the arithmetic Bohr radius within the framework of local Banach space theory, by establishing a connection between the asymptotic behaviors of these two radii. Among the various notions of Bohr radius, the arithmetic Bohr radius is more than a mere curious variant of the classical Bohr radius. It possesses rich structural properties; in particular, it plays an important role in describing the domain of existence of monomial expansions of bounded holomorphic functions on complete Reinhardt domains (see \cite{prengel-2005}). With this objective in mind, we first recall the following notion of the arithmetic Bohr radius, which is associated with the inequality \eqref{e-1.2}.
	\vspace{2mm}
	
	 Let $\Omega \subset \mathbb{C}^n$ be a given complete Reinhardt domain and $n \in \mathbb{N}$. Let $X$ be any complex Banach space. Let $\mathcal{F}(\Omega, X)$ be a set of $X$-valued holomorphic functions on $\Omega$. For any complex Banach spaces $X$ and $Y$, let $U:X\rightarrow Y$ be a bounded liner operator and $\norm{U} \leq \lambda$. For each $1\leq p<\infty$ and $\lambda\geq 1$, the $\lambda_{p}$-\textit{arithmetic Bohr radius} $A_{\lambda}(\mathcal{F}(\Omega, X), p, U)$ of $\Omega$ with respect to $\mathcal{F}(\Omega, X)$ is defined as \begin{equation*}  \sup \left\{\frac{1}{n}\sum_{i=1}^{n}r_i \,|\, r\in \mathbb{R}^{n}_{\geq 0},\, \mbox{for all}\, f\in \mathcal{F}(\Omega, X) : \sum_{\alpha} \norm{U(a_{\alpha})r^\alpha}^p_{Y} \leq \lambda^p\,\norm{f}^p_{\Omega,X}\right\}, 
	\end{equation*}
	where $\mathbb{R}^n_{\geq 0}=\{r=(r_1,.\,.\,.\,, r_n) \in \mathbb{R}^n: r_i\geq 0, 1\leq i\leq n\}.$ Let $H^{\infty}(\Omega,X)$ be the space of all bounded $X$-valued holomorphic functions on $\Omega$. We write $A_\lambda(\Omega,p,U)$ for $A_{\lambda}(H^{\infty}(\Omega,X),p,U)$ and $A(\Omega,p,U)$ for $A_1(\Omega,p,U)$. When $U=I:X\rightarrow X$ is the identity operator, we set $A_{\lambda}(\Omega, p,X):=A_{\lambda}(\Omega, p,I)$ and $A(\Omega, p,X):=A_{1}(\Omega, p,X)$. In the scalar valued case, we denote $A_{\lambda}(\Omega, p):=A_{\lambda}(\Omega, p,\mathbb{C})$ and $A(\Omega,p):=A_1(\Omega,p)$. 
The constant $A_{\lambda}(\Omega,1)$ was first considered by Defant {\it et al.} \cite{defant-2007} and is usually known as the classical arithmetic Bohr radius. Recently, the authors in \cite{	A-H-P-banach,Allu-Pal-arithnmetic} have studied the constant $A_\lambda(\Omega,p,U)$ and have investigated asymptotic estimates of it on the Minkowski spaces. 
\vspace{1mm}

In this paper, we investigate asymptotic estimates of the arithmetic Bohr radius on arbitrary bounded complete Reinhardt domains using techniques from local Banach space theory. In a recent work \cite{Himadri-local-Banach-1}, the author introduced a notion of Bohr radius for operator valued pluriharmonic functions defined on such domains and analyzed its asymptotic behavior using tools from local Banach space theory. Motivated by these developments, it is natural to ask whether the arithmetic Bohr radius can also be investigated in the setting of operator valued pluriharmonic functions. 

In order to address this problem, we begin by recalling the definitions of harmonic and pluriharmonic functions. A function $f:\Omega \rightarrow \mathbb{C}$ of class $C^2$ is called plurharmonic if its restriction to every complex line is harmonic. When $\Omega \subset\mathbb{C}^n$ is a simply connected domain in containing the origin, such a function is pluriharmonic precisely when it can be decomposed as $f=h+\overline{g}$, where $h,g$ are holomorphic on $\Omega$ with $g(0)=0$. This characterization admits a natural extension to the operator valued setting. In particular, if $\Omega$ is a simply connected complete Reinhardt domain, then a continuous function $f:\Omega \rightarrow \mathcal{B}(\mathcal{H})$ is pluriharmonic if and only it possesses a series representation of the form
\begin{equation} \label{e-1.3-a}
	f(z)=\sum_{m=0}^{\infty} \sum_{|\alpha|=m} a_{\alpha}\, z^{\alpha} + \sum_{m=1}^{\infty} \sum_{|\alpha|=m} b^{*}_{\alpha}\, \bar{z}^{\alpha},
\end{equation} 
where the associated functions $h(z)=\sum_{m=0}^{\infty} \sum_{\alpha} a_{\alpha}\, z^{\alpha}$ and $g(z)=\sum_{m=1}^{\infty} \sum_{\alpha} b_{\alpha}\, z^{\alpha}$ are $\mathcal{B}(\mathcal{H})$ valued holomorphic functions on $\Omega$, and $a_{\alpha}, b_{\alpha} \in \mathcal{B}(\mathcal{H})$. Here, $\mathcal{B}(\mathcal{H})$ denotes the Banach algebra of all bounded linear operators acting on a complex Hilbert space $\mathcal{H}$. For any $T \in \mathcal{B}(\mathcal{H})$, $\norm{T}$ will always denote the operator norm of $T$, and $T^{*}$ is the usual adjoint of $T$. Write $Re(T):=(T+T^*)/2$. Moreover, if $z=(z_1, \ldots, z_{n})$, then $\overline{z}$ denotes the componentwise complex conjugate $(\overline{z}_1, \ldots,\overline{z}_n)$. For such a domain $\Omega$ as above, we denote by $\mathcal{PH}(\Omega,X)$ the set of all bounded $X$-valued pluriharmonic functions on $\Omega$, where $X=\mathcal{B}(\mathcal{H})$. Boundedness is understood in the sense that $\norm{f}_{\Omega,X}:=\sup_{z \in \Omega}\,\norm{f(z)}_{X} < \infty$ for $f \in \mathcal{PH}(\Omega,X)$.
\vspace{2mm}

Let $Y$ be any complex Banach space. Let $\Omega\subset \mathbb{C}^n$ be a simply connected complete Reinhardt domain and $n\in \mathbb{N}$. Let $U:\mathcal{B}(\mathcal{H})\rightarrow Y$ be a bounded liner operator and $\norm{U} \leq \lambda$. For $1 \leq p < \infty$, the $\lambda$-powered Bohr radius of $U$, denoted by $R_{\lambda}(\Omega, p,U)$, is defined to be the supremum of all $r\geq 0$ such that for all $f \in \mathcal{PH}(\Omega,\mathcal{B}(\mathcal{H}))$ of the form \eqref{e-1.3-a} we have 
\begin{equation} \label{e-1.4-a}
	\sup_{z \in r\Omega}\,\sum_{m=0}^{\infty} \sum_{|\alpha|=m} (\norm{U(a_{\alpha})}^p_{Y} + \norm{U(b_{\alpha})}^p_{Y})|z^\alpha|^p \leq \lambda^p\,\norm{f}^p_{\Omega,\mathcal{B}(\mathcal{H})}. 
\end{equation}

\noindent Set $R(\Omega, p,U):=R_{1}(\Omega, p,U)$, $R_{\lambda}(\Omega, p,\mathcal{B}(\mathcal{H})):=R_{\lambda}(\Omega, p,U)$ whenever $U=I:\mathcal{B}(\mathcal{H})\rightarrow \mathcal{B}(\mathcal{H})$, $R(\Omega, p,\mathcal{B}(\mathcal{H})):=R_{1}(\Omega, p,\mathcal{B}(\mathcal{H}))$, $R_{\lambda}(\Omega, p):=R_{\lambda}(\Omega, p,\mathbb{C})$, and $R(\Omega,p):=R_1(\Omega,p)$. The Bohr phenomenon for harmonic functions in the scalar setting was first investigated in \cite{Abu-2010}. Subsequent developments extended this line of research to operator valued harmonic mappings (see \cite{bhowmik-2021}), and more recently to pluriharmonic functions taking values in Hilbert spaces (see \cite{hamada-Math-Nachr-2023}). The asymptotic behaviour of the constants $R(\Omega,1)$ and $R(\Omega,p)$ have been studied in \cite{hamada-JFA-2022,das-2024} whenever $\Omega=B_{\ell^n _q}$.
\vspace{2mm}

We now turn to the central theme of this article, namely the notion of the arithmetic Bohr radius for pluriharmonic functions. To the best of our knowledge, this concept has not been explored previously in the literature. This concept not been even studied in the simplest cases such as $U=I:\mathcal{B}(\mathcal{H})\rightarrow \mathcal{B}(\mathcal{H})$, $p=1$, $\lambda=1$, and for the classical domains such as the unit balls of Minkowski spaces, nor for the scalar valued harmonic and pluriharmonic functions. This gap motivates us to pose the following question. 
\begin{qsn} \label{qsn-1.2}
	Let $\Omega \subset \mathbb{C}^n$ be a simply connected complete Reinhardt domain. Can we study the arithmetic Bohr radius for operator valued pluriharmonic functions on $\Omega$?
\end{qsn}
We answer affirmatively to this question in this paper by introducing this notion in the following way.
Let $\mathcal{FP}(\Omega, \mathcal{B}(\mathcal{H}))$ be a set of $\mathcal{B}(\mathcal{H})$-valued pluriharmonic functions on $\Omega$. Let $U:\mathcal{B}(\mathcal{H})\rightarrow Y$ be a bounded liner operator and $\norm{U} \leq \lambda$.
	For each $1\leq p<\infty$ and $\lambda\geq 1$, the $\lambda_{p}$-\textit{arithmetic Bohr radius} $AP_{\lambda}(\mathcal{F}(\Omega, \mathcal{B}(\mathcal{H})), p, U)$ of $\Omega$ with respect to $\mathcal{FP}(\Omega, \mathcal{B}(\mathcal{H}))$ is defined as $\sup\left\{\frac{1}{n}\sum_{i=1}^{n}r_i:r\in \mathbb{R}^{n}_{\geq 0}\right\}$ such that
	\begin{equation*} 
\sum_{m=0}^{\infty} \sum_{|\alpha|=m} (\norm{U(a_{\alpha})}^p_{Y} + \norm{U(b_{\alpha})}^p_{Y})r^{p \alpha} \leq \lambda^p\norm{f}^p_{\Omega,\mathcal{B}(\mathcal{H})}
	\end{equation*}
for all $f\in \mathcal{FP}(\Omega, \mathcal{B}(\mathcal{H}))$ of the form \eqref{e-1.3-a}.
	Set $AP_\lambda(\Omega,p,U):=AP_\lambda(\mathcal{PH}(\Omega,\mathcal{B}(\mathcal{H})),p,U)$ and $AP(\Omega,P,U):=AP_1(\Omega,p,U)$. When $U=I:\mathcal{B}(\mathcal{H})\rightarrow \mathcal{B}(\mathcal{H})$ is the identity operator, we set $AP_{\lambda}(\Omega, p,\mathcal{B}(\mathcal{H})):=AP_{\lambda}(\Omega, p,I)$ and $AP(\Omega, p,X):=AP_{1}(\Omega, p,\mathcal{B}(\mathcal{H}))$. In the scalar valued case, we denote $AP_{\lambda}(\Omega, p):=AP_{\lambda}(\Omega, p,\mathbb{C})$ and $AP(\Omega,p):=AP_1(\Omega,p)$.  Let $AP^k_\lambda(\Omega,p,U):=AP_\lambda(\mathcal{PH}(^k\Omega,\mathcal{B}(\mathcal{H})),p,U)$. In the scalar valued case, we denote $AP^k_{\lambda}(\Omega, p):=AP^k_{\lambda}(\Omega, p,\mathbb{C})$ and $AP^k(\Omega,p):=A^k_1(\Omega,p)$.

	Recently, the present author in \cite{Himadri-local-Banach-1} has shown that the constant $R_{\lambda}(\Omega, p,U)$ is not necessarily nonzero for arbitrary operator $U$. Consequently, in the same paper, the author has established the condition when this constant is non-zero. Moreover, the author investigated asymptotic behaviors of this constant using invariants from local Banach space theory. With these developments in mind, it is natural to ask the following question for the constant $AP_{\lambda}(\Omega, p,U)$.
	\begin{qsn} \label{qsn-1.1}
		For a given simply connected complete Reinhardt domain $\Omega \subset \mathbb{C}^n$, can the the constant $AP_{\lambda}(\Omega, p,U)$ be studied in detail? If so, does it always remain nonzero, and what can be said about its asymptotic behavior?
	\end{qsn}
	One of our main contributions to this paper is to affirmatively answer the aforesaid question for arbitrary bounded simply connected complete Reinhardt domain by establishing the relation between Bohr radius and arithmetic Bohr radius. The following theorem relates the Bohr radius constant $R_{\lambda}(\Omega, p,U)$ and the arithmetic Bohr radius $AP_{\lambda}(\Omega, p,U)$ on any bounded  simply connected complete Reinhardt domain, which as a consequence shows the strictly positivity of $AP_{\lambda}(\Omega, p,U)$. This theorem actually help to obtain the lower estimate of $AP_{\lambda}(\Omega, p,U)$. A part of the proof follows a similar argument to that of \cite[Lemma 4.3]{defant-2012}; however, in that work the authors restricted attention only to scalar valued holomorphic functions.
	Before stating our result, we introduce the following standard notation. For bounded Reinhardt domains $\Omega_{1}, \Omega_{2}\subset \mathbb{C}^n$, let $S(\Omega_1, \Omega_{2}):= \inf \left\{s>0 : \Omega_{1} \subset s \, \Omega_{2}\right\}$. Recall that if $Z$ and $W$ are Banach sequence spaces then $S(B_{Z},B_{W})=\norm{Id:Z \rightarrow W}$.
	\begin{thm} \label{thm-1.8-atm}
			Let $Z=(\mathbb{C}^n, ||.||)$ be a Banach space such that $\chi(\{e_k\}^n_{k=1})=1$. Then 
		\begin{equation} \label{e-1.1-added}
			AP_\lambda(B_{Z},p,U) \geq \frac{\norm{Id:Z \rightarrow \ell^n_1}}{n} \, R_\lambda(B_{Z},p,U).
		\end{equation}
	Moreover, if $\norm{U}<\lambda$, then $AP_\lambda(B_{Z},p,U) \geq D.\frac{\norm{Id:Z \rightarrow \ell^n_1}}{n} \,. \frac{1}{\sup_{z \in B_{Z}}\norm{z}_{p}}$, where 
	$$
	D=\begin{cases}
		\max \left\{\frac{1}{4 \lambda\,2^{\frac{1}{p}}}\left(\frac{\lambda^p - \norm{U}^p}{2\lambda^p - \norm{U}^p}\right)^{\frac{1}{p}}\, , \frac{1}{4 \lambda\,2^{\frac{1}{p}}}\left(\frac{\lambda^p - \norm{U}^p}{\lambda^p - \norm{U}^p +1}\right)^{\frac{1}{p}}\, \frac{1}{\norm{U}}\right\}\,  & \text{for $\norm{U}\geq \frac{1}{4 \lambda\,2^{\frac{1}{p}}}$},\\[3mm]
		\max \left\{\frac{1}{4 \lambda\,2^{\frac{1}{p}}}\left(\frac{\lambda^p - \norm{U}^p}{2\lambda^p - \norm{U}^p}\right)^{\frac{1}{p}}, \frac{1}{4 \lambda\,2^{\frac{1}{p}}}\left(\frac{\lambda^p - \norm{U}^p}{\lambda^p - \norm{U}^p +1}\right)^{\frac{1}{p}}\right\} & \text{for $0<\norm{U}< \frac{1}{4 \lambda\,2^{\frac{1}{p}}}$}.
	\end{cases}
	$$
	\end{thm}
Observing that $S(\Omega,B_{\ell^n_1})=\sup_{z \in \Omega} \norm{z}_{1}$, and following arguments similar to those used in the proof of Theorem~\ref{thm-1.8-atm}, it is straightforward to verify that the same conclusion holds for any simply connected complete Reinhardt domain. We therefore state the result and omit the details to avoid repetition.
\begin{cor} \label{cor-1.1-added}
Let $\Omega$ be a simply connected complete Reinhardt domain in $\mathbb{C}^n$. Then 
\begin{equation*}
	AP_\lambda(\Omega,p,U) \geq \frac{S(\Omega,B_{\ell^n_1})}{n} \, R_\lambda(\Omega,p,U).
\end{equation*}
Moreover, if $\norm{U}<\lambda$, then $AP_\lambda(\Omega,p,U) \geq D.\frac{S(\Omega,B_{\ell^n_1})}{n} \,. \frac{1}{\sup_{z \in \Omega}\norm{z}_{p}}$, where $D$ as in Theorem \ref{thm-1.8-atm}.
\end{cor}
Applying this corollary to the case $\Omega=B_{\ell^n_q}$, the Minkowski space, yields the following result. 
\begin{cor}
Let $Y$ be any complex Banach space and $U:\mathcal{B}(\mathcal{H})\rightarrow Y$ a bounded liner operator with $\norm{U} \leq \lambda$. For each $1\leq p<\infty$ and $1\leq q \leq\infty$,
\begin{equation*}
AP_\lambda(B_{\ell^n_q},p,U) \geq D\,
\begin{cases}	
	n^{- \frac{1}{q}} & \text{for $q\leq p$}, \\[2mm]
	 n^{ - \frac{1}{p}} & \text{for $q>p$},
\end{cases}
\end{equation*}
where $D$ as in Theorem \ref{thm-1.8-atm}.
\end{cor}
	In particular, we obtain the following lower bound of the $\lambda$-powered arithmetic Bohr radius for operator-valued pluriharmonic functions whenever $U$ is the identity operator on $X=\mathcal{B}(\mathcal{H})$.
	\begin{cor}
		Let $X=\mathcal{B}(\mathcal{H})$ and $1<\lambda$. Then for all $n\in \mathbb{N}$ and $1\leq p < \infty$, we have
		\begin{equation*}
			AP_{\lambda}(\Omega, p,X) \geq \frac{1}{2^{2+\frac{1}{p}}} \frac{\left(\lambda^p -1\right)^{\frac{1}{p}}}{\lambda^2}\, \frac{S(\Omega, B_{\ell^n_1})}{n}\, \frac{1}{\sup_{z \in \Omega}\norm{z}_{p}}.
		\end{equation*}
	\end{cor}
	\begin{rem}
	In view of Theorem \ref{thm-1.8-atm} and Corollary \ref{cor-1.1-added}, it is evident that the constant $AP_\lambda(\Omega,p,U)$ is non-zero whenever $\norm{U}<\lambda$. Moreover, a direct comparison of the definitions of arithmetic Bohr radius and Bohr radius shows that $R_\lambda(\Omega,p,U)\leq AP_\lambda(\Omega,p,U)$. As a consequence, the vanishing of $AP_\lambda(\Omega,p,U)$ necessarily implies the vanishing of $R_\lambda(\Omega,p,U)$. However, the converse implication fails in general, reflecting a genuine distinction between these two notions. In the recent work, the author in \cite{Himadri-local-Banach-1} has shown that $R_\lambda(\Omega,p,U)=0$ when $\norm{U}=\lambda$. This result naturally leads to the question of whether the arithmetic Bohr radius exhibits the same behavior, namely, whether $AP_\lambda(\Omega,p,U)$ also vanishes when $\norm{U}=\lambda$. The example presented below shows that this is not the case and the constant $AP_{\lambda}(\Omega, p,U)$ may be strictly positive. 
		For a bounded simply connected complete Reinhardt domain $\Omega \subseteq \mathbb{C}^n$ and $k \in \mathbb{N}$, consider the functions  
		$F_{k}: \Omega \rightarrow \mathcal{B}(\mathcal{H})$ by $$F_{k}(z)=(i\, \cos \frac{1}{k}) I_{\mathcal{H}} + (\frac{1}{2} \sin \frac{1}{k})  I_{\mathcal{H}}z_{1}+(\frac{1}{2} \sin \frac{1}{k})  I_{\mathcal{H}} \overline{z_{1}}$$
		for $z=(z_{1}, \ldots,z_{n}) \in \Omega$, where $I_{\mathcal{H}}$ is the identity on $\mathcal{H}$. Let $U=\lambda \, I$, with $I$ the identity on $\mathcal{B}(\mathcal{H})$. Then $F_{k}(0)=(i\, \cos \frac{1}{k})  I_{\mathcal{H}}$. Clearly, $F_{k} \in \mathcal{PH}(\Omega,\mathcal{B}(\mathcal{H}))$ and $\norm{F_{k}}_{\Omega,\mathcal{B}(\mathcal{H})} \leq 1$. Let if possible, there exists $r_{0}>0$ such that \eqref{e-1.4-a} holds for all $f \in \mathcal{PH}(\Omega,\mathcal{B}(\mathcal{H}))$ for all $z \in r_{0}\, \Omega$. Applying this to the functions $F_{k}$ yields, 
		\begin{equation} \label{e-example-inequality}
		\lambda^p \, |\cos \frac{1}{k}|^p + \lambda^p \, |\sin \frac{1}{k}|^p\, |z_{1}|^p \leq \lambda^p
		\end{equation}
		all $z \in r_{0}\, \Omega$ and for all $k \in \mathbb{N}$. However, this inequality fails for sufficiently large $k$, since $\cos (1/k) \rightarrow 1$ and $\sin(1/k) \rightarrow 0$ as $k \rightarrow\infty$, while $z$
		can be chosen in $r_{0}\, \Omega$ such that $z_{1} \neq 0$. This fact shows that there is no $r_{0}>0$ such that \eqref{e-example-inequality} holds, and hence $R_\lambda(\Omega,p,U)=0$. On the other hand, it is easy to check that \eqref{e-example-inequality} holds for all the the points $(0,\xi_2, \ldots, \xi_n) \in \Omega$ such that $\xi_j \geq 0$, $j=1,\ldots,n$. Thus, by the definition, $AP_\lambda(\Omega,p,U)>0$. Hence, we conclude that $R_\lambda(\Omega,p,U)=0$ while $AP_\lambda(\Omega,p,U)>0$ whenever $U=\lambda\, I$, thereby highlighting a fundamental difference between the arithmetic Bohr radius and the Bohr radius.
	\end{rem}
The following result links between the arithmetic Bohr radii of the spaces $\mathcal{PH}(\Omega,\mathcal{B}(\mathcal{H}))$ and $\cup^{\infty}_{m=1} \mathcal{PH}(^m \Omega,\mathcal{B}(\mathcal{H}))$. More precisely, we estimate $AP_{\lambda}(\Omega, p,U)$ in terms of $AP_{\lambda}\big(\cup^{\infty}_{m=1} \mathcal{PH}(^m \Omega,X), p,U\big)$. 
\begin{thm} \label{tthm-1.8-atm}
	Let $\Omega$ be a simply connected complete Reinhardt domain in $\mathbb{C}^n$. Let $X=\mathcal{B}(\mathcal{H})$ and $Y$ any complex Banach space, and $U: X \rightarrow Y$ be a non-null bounded linear operator such that $\norm{U}< \lambda$. Then, for all $p\in[1,\infty)$, $\lambda>1$, and $n \in \mathbb{N}$, we have\\
	$(1)$ $\left(\frac{\lambda^p - \norm{U}^p}{2\lambda^p - \norm{U}^p}\right)^{\frac{1}{p}}\, AP_{\lambda}\big(\cup^{\infty}_{m=1} \mathcal{PH}(^m \Omega,X), p,U\big) \leq AP_{\lambda}(\Omega, p,U) \leq AP_{\lambda}\big(\cup^{\infty}_{m=1} \mathcal{PH}(^m \Omega,X), p,U\big)$ \\
	$(2)$ $\left(\frac{\lambda^p - \norm{U}^p}{\lambda^p - \norm{U}^p +1}\right)^{\frac{1}{p}}\, AP\big(\cup^{\infty}_{m=1} \mathcal{PH}(^m \Omega,X), p,U\big) \leq AP_{\lambda}(\Omega, p,U) \leq \lambda \, \inf_{m\in \mathbb{N}} AP\big( \mathcal{PH}(^m \Omega,X), p,U\big)$.
\end{thm}

We now impose additional conditions on the initial coefficient, which lead to improve estimates for the arithmetic Bohr radius obtained in Theorem \ref{tthm-1.8-atm}. More precisely, we consider the following class of functions. Let $\mathcal{PH}^{*}(\Omega,\mathcal{B}(\mathcal{H})):=\{f \in \mathcal{PH}(\Omega,\mathcal{B}(\mathcal{H})): f(0)=\mu\, I, \, \mu \in [0, \norm{f}_{\Omega,\mathcal{B}(\mathcal{H})}]\}$. Denote $AP^{*}_{\lambda}(\Omega, p,U):= AP_{\lambda}(\mathcal{PH}^{*}(\Omega,\mathcal{B}(\mathcal{H})), p,U)$.
\begin{thm} \label{tthm-1.9-atm}
	Let $X=\mathcal{B}(\mathcal{H})$ and $Y$ any complex Banach space, and $U: X \rightarrow Y$ be a non-null bounded linear operator such that $\norm{U}< \lambda$. Then, for all $p\in[1,\infty)$, $\lambda>1$, and $n \in \mathbb{N}$, we have\\
	$(1)$ $C\, AP_{\lambda}\big(\cup^{\infty}_{m=1} \mathcal{PH}(^m \Omega,X), p,U\big) \leq AP^{*}_{\lambda}(\Omega, p,U) \leq AP_{\lambda}\big(\cup^{\infty}_{m=1} \mathcal{PH}(^m \Omega,X), p,U\big)$ \\
	$(2)$ $D\, AP\big(\cup^{\infty}_{m=1} \mathcal{PH}(^m \Omega,X), p,U\big) \leq AP^{*}_{\lambda}(\Omega, p,U) \leq \lambda \, \inf_{m\in \mathbb{N}} AP\big( \mathcal{PH}(^m \Omega,X), p,U\big)$.
	Here,
	$
	C=
	\max \left\{\frac{1}{(1+8^p)^{1/p}}, \left(\frac{\lambda^p - \norm{U}^p}{2\lambda^p - \norm{U}}\right)^{\frac{1}{p}}\right\}
	$
	and
	$D=
	\max \left\{\frac{\lambda}{(\lambda^p+8^p)^{1/p}}, \left(\frac{\lambda^p - \norm{U}^p}{\lambda^p - \norm{U}^p +1}\right)^{\frac{1}{p}}\right\}.
	$
\end{thm}
We now turn to the second part of Question \ref{qsn-1.1}, which is the central focus of this paper, namely the asymptotic estimate of the constant $AP_{\lambda}(\Omega, p,U)$ through invariants from local Banach space theory, such as unconditional bases and projection constants. 

In order to study this, we shall make use of the following correspondence.  Recall that a Schauder basis $\{w_{k}\}$ of a Banach space $W$ is said to be $1$-unconditional if $\chi(\{w_k\})=1$ (see more details in Section $3$). We usually denote the canonical basis vectors of $Z=(\mathbb{C}^n, ||.||)$ by $e_{k}$, $k=1, \ldots,n$. 
  Let $Z=(\mathbb{C}^n, ||.||)$ be an $n$-dimensional Banach space for which the canonical basis vectors $e_{k}$ form a normalized $1$-unconditional basis (or equivalently, the open unit ball $B_{Z}$ is a complete Reinhardt domain). There is a one-to-one correspondence between bounded convex complete Reinhardt domains $\Omega \subseteq\mathbb{C}^n$ and the open unit balls of norms on $\mathbb{C}^n$ for which the canonical basis vectors $e_{k}$ form a normalized $1$-unconditional basis. Indeed, the Minkowski functional $p_{\Omega}:\mathbb{C}^n \rightarrow \mathbb{R}_{+}$ of $\Omega$ defined by $$p_{\Omega}(z):= \inf \left\{t>0 : z/t \in \Omega\right\}$$
	is a norm on $\mathbb{C}^n$, and $\Omega$ coincides with the open unit ball of $(\mathbb{C}^n,p_{\Omega})$. Conversely, if $Z=(\mathbb{C}^n, ||.||)$ is a Banach space whose canonical basis is normalized $1$-unconditional, then its open unit ball $B_{	Z}$ is clearly a bounded convex complete Reinhardt domain in $\mathbb{C}^n$. Thus, the study of the constant $R_{\lambda}(\Omega, p,U)$ for the unit balls $\Omega=B_{Z}$ of finite dimensional complex Banach spaces $(\mathbb{C}^n, ||.||)$ with the normalized $1$-unconditional canonical bases, is equivalent to studying $R_{\lambda}(\Omega, p,U)$ over bounded convex complete Reinhardt domains $\Omega \subseteq\mathbb{C}^n$.
	
	Our investigation naturally separates into two cases according to whether $\mathcal{B}(\mathcal{H})$ is finite and infinite dimensional. In each case, taking into account the above correspondence, we estimate the asymptotic behavior of $AP_{\lambda}(\Omega, p,U)$ in the following steps. First, we establish the result for $\Omega=B_{Z}$, that is, for a convex complete Reinhardt domain. Next, by combining this estimate with the above correspondence and Lemma \ref{lem-3.6}, we derive the asymptotic behavior of $AP_{\lambda}(\Omega, p,U)$ for every bounded simply connected complete Reinhardt domain $\Omega \subseteq \mathbb{C}^n$, not necessarily convex.
	
	 We start with the finite-dimensional case. The following theorem determines the exact asymptotic behavior of the arithmetic Bohr radius $AP_\lambda(B_{Z},p,\mathcal{B}(\mathcal{H}))$ when $U$ is the identity operator on $\mathcal{B}(\mathcal{H})$.
	\begin{thm} \label{thm-1.2-atm}
		Let $Z=(\mathbb{C}^n, ||.||)$ be a Banach space such that $\chi(\{e_k\}^n_{k=1})=1$. Let $X=\mathcal{B}(\mathcal{H})$ be a finite dimensional complex Banach space and $\lambda>1$. Then 
		\comment{\begin{equation}
				R^{m}_{\lambda}(B_{Z}, p) \geq	\begin{cases}
					E_{1}(X)\, \lambda^{\frac{1}{m}} \, \max \left\{\frac{1}{\sqrt{n}\, \norm{Id: \ell^n_{2} \rightarrow Z}}, \, \frac{1}{\left(\frac{m^m}{m!}\right)^{1/m}\, \norm{Id: Z \rightarrow \ell^n_{1}}}\right\} & \text{for $p=1$}, \\[3mm]
					2^{-\frac{p+5}{pm}}\,  \lambda^{\frac{1}{m}}\, \norm{I:\ell^n _2 \rightarrow Z}^{-\frac{2}{p}} \,\norm{I-\real (a_{0})}^{-\frac{2}{pm}} & \text{for $p \geq 2$}, \\[2mm]
					
					E_{2}(X)\, \lambda^{\frac{1}{m}} \,\max \left\{\frac{1}{(\sqrt{n})^{1-\theta}\, \norm{Id: \ell^n_{2} \rightarrow Z}}, \, \frac{\norm{I:\ell^n _2 \rightarrow Z}^{-\theta}}{\left(\frac{m^m}{m!}\right)^{(1-\theta)/m}\, \norm{Id: Z \rightarrow \ell^n_{1}}^{1-\theta}}\right\} & \text{for $1<p< 2$}, 
				\end{cases}
			\end{equation}
			where $\theta=(2(p-1))/p$. Moreover,}
		\begin{equation*}
			AP_{\lambda}(B_{Z}, p,X) \geq	\begin{cases}
				E_{1}(X)\,  \max \left\{\frac{\norm{Id: Z \rightarrow \ell^n_{1}}}{n^{3/2}\, \norm{Id: \ell^n_{2} \rightarrow Z}}, \, \frac{1}{e\,n}\right\}\, \left(\frac{\lambda^p -1}{2\lambda^p - 1}\right)^{\frac{1}{p}} & \text{for $p=1$}, \\[2mm]
				E_{3} \frac{\norm{Id: Z \rightarrow \ell^n_{1}}}{n}\norm{I:\ell^n _2 \rightarrow Z}^{-\frac{2}{p}}\, \left(\frac{\lambda^p -1}{2\lambda^p - 1}\right)^{\frac{1}{p}}\, & \text{for $p \geq 2$}, \\[2mm]
				
				E_{2}(X)\,  \max \left\{\frac{\norm{Id: Z \rightarrow \ell^n_{1}}}{n^{\frac{3-\theta}{2}} \norm{Id: \ell^n_{2} \rightarrow Z}},  \frac{\norm{I:\ell^n _2 \rightarrow Z}^{-\theta}}{e^{(1-\theta)}n \norm{Id: Z \rightarrow \ell^n_{1}}^{-\theta}}\right\}\left(\frac{\lambda^p -1}{2\lambda^p - 1}\right)^{\frac{1}{p}} \hspace{-2mm} & \text{for $1<p< 2$}, 
			\end{cases}
		\end{equation*}
		where $\theta = (2(p-1))/p$. Here, $E_{1}(X)$, $E_{2}(X)$, and $E_{3}$ are positive constants: $E_{1}(X)$ depends only on $X$, $E_{2}(X)$ depends on both $X$ and $p$, while $E_{3}$ is independent of $X$ and depends only on $p$. Furthermore,
		for any $1\leq q \leq \infty$,
		\begin{equation*}
			AP_{\lambda}(B_{Z}, p,X) \leq d\,\, \lambda^{\frac{1}{\log\,n}}\, \norm{Id:Z \rightarrow \ell^n_q}\, (\log\,n)^{1-\frac{1}{\min\{2,q\}}}\,\,\left(\frac{1}{n}\right)^{\frac{1}{p}+\frac{1}{\max\{2,q\}}-\frac{1}{2}}
		\end{equation*}
		for some constant $d>0$.
	\end{thm}
In particular, for the Minkowski spaces $Z=\ell^n _q$, we obtain the following lower bound. The upper bound is obtained in the proof of Theorem \ref{thm-1.2-atm}, following the technique from \cite{defant-2007}.
\begin{cor} \label{cor-1.4}
		Let $X=\mathcal{B}(\mathcal{H})$ be finite dimensional and $\lambda>1$. Then, for $1 \leq q \leq \infty$,
	\begin{equation*}
		AP_{\lambda}(B_{\ell^n _q}, p,X) \geq	\begin{cases}
			E'_{3}(X)\,\left(\frac{\lambda^p -1}{2\lambda^p - 1}\right)^{\frac{1}{p}}\, n^{-\frac{1}{q}} \, \left(\frac{\log n}{n}\right)^{ \left(1- \frac{1}{\min\{q,2\}}\right)} & \text{for $p=1$}, \\[2mm]
			E'_{2} \left(\frac{\lambda^p -1}{2\lambda^p - 1}\right)^{\frac{1}{p}}\, n^{-\left(\frac{1}{p}+\frac{1}{q}\right)} & \text{for $p \geq q$}, \\[2mm]
			
			E_{4}(X)\,  \left(\frac{\lambda^p -1}{2\lambda^p - 1}\right)^{\frac{1}{p}} \,n^{-\frac{p(q-1)+q(p-1)}{pq(q-1)}}\, \left( \frac{\log n}{n}\right)^{\left(1- \frac{1}{\min\{q,2\}}\right)\, \frac{q-p}{p(q-1)}} & \text{for $1<p< q$}. 
		\end{cases}
	\end{equation*}
	Here, $E'_{3}(X)$, $E'_{2}$, and $E_{4}(X)$ are positive constants: $E'_{3}(X)$ depends only on $X$, $E_{4}(X)$ depends on both $X$ and $p$, while $E'_{2}$ is independent of $X$ and depends only on $p$
\end{cor}
Theorem \ref{thm-1.2-atm} leads to the following asymptotic estimates of the arithmetic Bohr radius for Banach sequence spaces, whose broad applications are described in Section $2$.
\begin{thm} \label{thm-1.3-a}
	Let $Z$ be a Banach sequence space, $X=\mathcal{B}(\mathcal{H})$ be finite dimensional, and $\lambda>1$. For $n\in \mathbb{N}$, let $Z_{n}$ be the linear span of $\{e_{k}, \, k=1,\ldots,n\}$. Then
	\begin{enumerate}
		\item if $Z$ is a subset of $\ell_2$ we have \begin{equation*}
			AP_{\lambda}(B_{Z}, p,X) \geq	\begin{cases}
				E_{1}(X)\,  \max \left\{\frac{\norm{\sum_{k=1}^{n}e^{*}_{k}}_{Z^{*}}}{n^{3/2}\, \norm{Id: \ell^n_{2} \rightarrow Z}}, \, \frac{1}{e\,n}\right\}\, \left(\frac{\lambda^p -1}{2\lambda^p - 1}\right)^{\frac{1}{p}} & \text{for $p=1$}, \\[2mm]
				E_{3} \frac{\norm{\sum_{k=1}^{n}e^{*}_{k}}_{Z^{*}}}{n}\norm{Id:\ell^n _2 \rightarrow Z}^{-\frac{2}{p}}\, \left(\frac{\lambda^p -1}{2\lambda^p - 1}\right)^{\frac{1}{p}}\, & \text{for $p \geq 2$}, \\[2mm]
				
				E_{2}(X)\,  \max \left\{\frac{\norm{\sum_{k=1}^{n}e^{*}_{k}}_{Z^{*}}}{n^{\frac{3-\theta}{2}} \norm{Id: \ell^n_{2} \rightarrow Z}},  \frac{\norm{Id:\ell^n _2 \rightarrow Z}^{-\theta}}{e^{(1-\theta)}n \norm{\sum_{k=1}^{n}e^{*}_{k}}_{Z^{*}}^{-\theta}}\right\}\left(\frac{\lambda^p -1}{2\lambda^p - 1}\right)^{\frac{1}{p}} \hspace{-2mm} & \text{for $1<p< 2$}, 
			\end{cases}
		\end{equation*}
		and 
		$$AP_{\lambda}(B_{Z_{n}}, p,X) \leq d\,  \lambda^{\frac{1}{\log\,n}}\,\, \sqrt{\log \, n} \, \left(\frac{1}{n}\right)^{\frac{1}{p}};$$
		\item if $Z$ is symmetric and $2$-convex, we have
		\begin{equation*}
			AP_{\lambda}(B_{Z}, p,X) \geq	\begin{cases}
				E_{1}(X)\,  \max \left\{\frac{\norm{\sum_{k=1}^{n}e^{*}_{k}}_{Z^{*}}}{n^{3/2}}, \, \frac{1}{e\,n}\right\}\, \left(\frac{\lambda^p -1}{2\lambda^p - 1}\right)^{\frac{1}{p}} & \text{for $p=1$}, \\[2mm]
				E_{3} \frac{\norm{\sum_{k=1}^{n}e^{*}_{k}}_{Z^{*}}}{n}\, \left(\frac{\lambda^p -1}{2\lambda^p - 1}\right)^{\frac{1}{p}}\, & \text{for $p \geq 2$}, \\[2mm]
				
				E_{2}(X)\,  \max \left\{\frac{\norm{\sum_{k=1}^{n}e^{*}_{k}}_{Z^{*}}}{n^{\frac{3-\theta}{2}}},  \frac{1}{e^{(1-\theta)}n \norm{\sum_{k=1}^{n}e^{*}_{k}}_{Z^{*}}^{-\theta}}\right\}\left(\frac{\lambda^p -1}{2\lambda^p - 1}\right)^{\frac{1}{p}} \hspace{-2mm} & \text{for $1<p< 2$}, 
			\end{cases}
		\end{equation*}
		and 
		$$AP_{\lambda}(B_{Z}, p,X) \leq d\,  \lambda^{\frac{1}{\log\,n}}\, \norm{\sum_{k=1}^{n}e^{*}_{k}}_{Z^{*}}\,\,\, \sqrt{\log \, n}\,\, \left(\frac{1}{n}\right)^{\frac{1}{p}+\frac{1}{2}}.
		$$
	\end{enumerate}
	The constants $E_{1}, E_{2}, E_{3}$, and $d$ are same as in Theorem \ref{thm-1.2-atm}.
\end{thm}
To determine the exact asymptotic behavior of the constant $AP_{\lambda}(B_{Z}, p,\mathcal{B}(\mathcal{H}))$ in the case where $\mathcal{B}(\mathcal{H})$ is infinite dimensional, we make essential use of Lemma \ref{lem-3.6}. This lemma allows us to reduce the problem to the study of this constant for such a domain in which the estimation of bounds becomes significantly more simpler. Among complete Reinhardt domains in $\mathbb{C}^n$, the most natural and widely studied examples are the unit balls $B_{\ell^n _q}$ of the Minkowski spaces $\ell^n _q$, $1\leq q \leq \infty$. Consequently, we begin by analyzing the behavior of $AP_{\lambda}(\Omega, p,\mathcal{B}(\mathcal{H}))$ when $\Omega=B_{\ell^n _q}$. 

In the finite-dimensional setting, Corollary \ref{cor-1.4} shows that the asymptotics of this constant necessarily involve a logarithmic factor. By contrast, the situation changes fundamentally in the infinite-dimensional case. As we show in the theorem below, the logarithmic term disappears entirely, and the asymptotic decay of $AP_{\lambda}(B_{\ell^n _q}, p,\mathcal{B}(\mathcal{H}))$ is governed instead by the geometric properties of the Banach space $X$, namely its optimal cotype $Cot(X)$(see Section $2$ for the definition), where $X=\mathcal{B}(\mathcal{H})$. 
	We now obtain the asymptotic estimates of the arithmetic Bohr radius for $B_{\ell^n _q}$ whenever $X=\mathcal{B}(\mathcal{H})$ is infinite dimensional.	
	\begin{thm} \label{thm-1.3-atm}
		Let $X=\mathcal{B}(\mathcal{H})$ be an infinite dimensional complex Banach space  and $\lambda>1$. Then $AP_{\lambda}(B_{\ell^n _q}, p,X) \geq \tilde{\Psi}_{1}(n,\lambda,p,q)$, where
		\begin{equation*}
			\tilde{\Psi}_{1}(n,\lambda,p,q):=	\begin{cases}
				
				E_{5}\,\dfrac{\left(\lambda^p -1\right)^{\frac{1}{p}}}{\lambda} \frac{1}{n^{1/q}} & \text{for $p>\min \{Cot(X),q\}$}, \\[2mm]
				E_{6}\, \dfrac{\left(\lambda^p -1\right)^{\frac{1}{p}}}{\lambda}\, n^{{\frac{1}{\min\{t,q\}}}-\frac{1}{p}-\frac{1}{q}} & \text{for $p\leq \min \{Cot(X),q\}$},
			\end{cases}
		\end{equation*}
		if $X$ has a finite cotype $t$. Here $E_{5}>0$ and $E_{6}>0$ are constant independent of $X$ and $n$. 
		On the other hand, we have
		\begin{equation*}
			AP_{\lambda}(B_{\ell^n _q}, p,X) \leq \tilde{\Psi}_{2}(n,\lambda,p,q):=	\begin{cases}	
				2^{1-\frac{1}{p}}\, \lambda\, n^{- \frac{1}{p}} & \text{for $q\leq Cot(X)$}, \\[2mm]
				2^{1-\frac{1}{p}}\,\lambda\, n^{\frac{1}{Cot(X)}\, - \,\frac{1}{p}-\frac{1}{q}} & \text{for $q>Cot(X)$}.
			\end{cases}
		\end{equation*}
	\end{thm}
	
	We are now in a stage to give the bounds of the arithmetic Bohr radius $AP_{\lambda}(B_{Z}, p,X)$ for the unit ball $B_{Z}$ of $Z=(\mathbb{C}^n,||.||)$ whenever $X=\mathcal{B}(\mathcal{H})$ is a infinite dimensional complex Banach space. In view of Theorem \ref{thm-1.3-atm} and Lemma \ref{lem-3.6}, we obtain the following result.
	\begin{thm} \label{thm-1.4-atm}
		Let $Z=(\mathbb{C}^n, ||.||)$ be a Banach space such that $\chi(\{e_k\}^n_{k=1})=1$. Let $\mathcal{B}(\mathcal{H})$ be infinite dimensional, and $\tilde{\Psi}_{1}(n,\lambda,p,q)$ and $\tilde{\Psi}_{2}(n,\lambda,p,q)$ be as in Theorem \ref{thm-1.3-atm}. For $\lambda>1$, 
		$$\frac{\tilde{\Psi}_{1}(n,\lambda,p,q)}{\norm{Id: \ell^n_{q} \rightarrow Z}} \, \leq AP_{\lambda}(B_{Z}, p,X) \leq \norm{Id: Z \rightarrow \ell^n_{q}} \tilde{\Psi}_{2}(n,\lambda,p,q).
		$$
	\end{thm}
The correspondence preceding Theorem \ref{thm-1.2-atm} allows Theorem \ref{thm-1.4-atm} to be extended to bounded simply connected convex complete Reinhardt domains. Moreover, combining Theorem \ref{thm-1.3-atm} with Lemma \ref{lem-3.6} yields the asymptotic estimates for all bounded simply connected complete Reinhardt domains, not necessarily convex. We record this result as a corollary.
\begin{cor} \label{cor-1.4-atm}
	Let $\Omega$ be a bounded simply connected complete Reinhardt domain. Let $\mathcal{B}(\mathcal{H})$ be infinite dimensional, and $\tilde{\Psi}_{1}(n,\lambda,p,q)$ and $\tilde{\Psi}_{2}(n,\lambda,p,q)$ be as in Theorem \ref{thm-1.3-atm}. For $\lambda>1$, 
	$$\frac{\tilde{\Psi}_{1}(n,\lambda,p,q)}{S(B_{\ell^n _q},\Omega)} \, \leq AP_{\lambda}(\Omega, p,X) \leq S(\Omega,B_{\ell^n _q}) \tilde{\Psi}_{2}(n,\lambda,p,q).
	$$
\end{cor}
An application of Theorem \ref{thm-1.4-atm}, together with the norm estimates for the identity operators mentioned in the proof of Theorem \ref{thm-1.3-a}, yields the following result.
	\begin{thm} \label{thm-1.6}
		Let $Z$ be a Banach sequence space and $X=\mathcal{B}(\mathcal{H})$ be infinite dimensional. For $n\in \mathbb{N}$, let $Z_{n}$ be the linear span of $\{e_{k}, \, k=1,\ldots,n\}$. Then for $\lambda>1$
		\begin{enumerate}
			\item if $Z$ is subset of $\ell_2$, we have 
			$$\frac{\tilde{\Psi}_{1}(n,\lambda,p,2)}{\norm{Id: \ell^n_{2} \rightarrow Z_n}} \, \leq AP_{\lambda}(B_{Z_n}, p,X) \leq  \tilde{\Psi}_{2}(n,\lambda,p,2);
			$$
			\item if $Z$ is symmetric and $2$-convex, we have $$\tilde{\Psi}_{1}(n,\lambda,p,2)\, \leq AP_{\lambda}(B_{Z}, p,X) \leq  \frac{\norm{\sum_{m=1}^{n}e^{*}_{m}}_{Z^{*}}}{\sqrt{n}}\, \tilde{\Psi}_{2}(n,\lambda,p,2).
			$$
		\end{enumerate}
	\end{thm}
	We have so far investigated the asymptotic behavior of the arithmetic Bohr radius associated with the identity operator $U:\mathcal{B}(\mathcal{H}) \rightarrow \mathcal{B}(\mathcal{H})$. We now turn to the study of this constant for a general operator $U$, which need not be the identity. To this end, we first analyze the Bohr radius constant corresponding to an arbitrary operator. Using this analysis, we then derive bounds for the arithmetic Bohr radius.
	\vspace{1mm}
	
	 With this objective in mind, we begin by establishing the following result, which may be viewed as an analogue of \cite[Theorem 5.3]{defant-2012}. While the authors there proved the result for the unit polydisc, our approach extends it to the unit ball of an arbitrary Banach space admitting a $1$-unconditional basis.
	\begin{thm} \label{thm-1.8}
	Let $Y$ be a $s$-concave Banach lattice, with $2\leq s < \infty$, and $U: X \rightarrow Y$ an $(t,1)$-summing operator with $1 \leq t \leq s$. Then there is a constant $C>0$ such that
	\begin{equation*}
	\left(\sum_{|\alpha|=m} \norm{U(a_{\alpha})z^{\alpha}}^{\beta}_{Y}\right)^{\frac{1}{\beta}} \leq C^m \, \norm{P}_{B_{Z},X}
	\end{equation*}
holds for every homogeneous polynomial $P(z)=\sum_{|\alpha|=m} a_{\alpha}\, z^{\alpha}$ on $Z$, where 
\begin{equation*}
\beta= \frac{mst}{s+ (m-1)t}.
\end{equation*}
	\end{thm}

We are now in a position to establish lower bounds for both the Bohr radius and the arithmetic Bohr radius for a general operator.
	\begin{thm} \label{thm-1.9}
	Let $1 \leq p < \infty$ and $1\leq q \leq \infty$.
	\begin{enumerate}
		\item Let $\mathcal{H}_1$ and $\mathcal{H}_2$ be two complex Hilbert spaces, and $U: \mathcal{B}(\mathcal{H}_1) \rightarrow \mathcal{B}(\mathcal{H}_2)$ be a non-null bounded linear operator with $\norm{U}<\lambda$. If $\mathcal{B}(\mathcal{H}_1)$ or $\mathcal{B}(\mathcal{H}_2)$ is of cotype $t$ with $2 \leq t \leq \infty$, then for every $n \in \mathbb{N}$,	
	\begin{align*}
	R_{\lambda}(B_{Z}, p,U) \geq \frac{\Phi_{1}(n,\lambda,p,q)}{\norm{Id:Z \rightarrow \ell^n _q} \norm{Id:\ell^n _q \rightarrow Z}}
	\end{align*} 
and
\begin{equation*}
	AP_{\lambda}(B_{Z}, p,U) \geq \frac{\Phi_{1}(n,\lambda,p,q)\, \norm{Id:Z \rightarrow \ell^n _1}}{n\norm{Id:Z \rightarrow \ell^n _q} \norm{Id:\ell^n _q \rightarrow Z}},
\end{equation*}
	where
	\begin{equation*}
		\Phi_{1}(n,\lambda,p,q):=	\begin{cases}
			
			E_{5}\,\dfrac{\left(\lambda^p -\norm{U}^p\right)^{\frac{1}{p}}}{\lambda} & \text{for $p>\min \{Cot(\mathcal{B}(\mathcal{H}_1)),Cot(\mathcal{B}(\mathcal{H}_2)),q\}$}, \\[2mm]
			E_{6}\, \dfrac{\left(\lambda^p -\norm{U}^p\right)^{\frac{1}{p}}}{\lambda}\,\, n^{{\frac{1}{\min\{t,q\}}}-\frac{1}{p}} & \text{for $p\leq \min \{Cot(\mathcal{B}(\mathcal{H}_1)),Cot(\mathcal{B}(\mathcal{H}_2)),q\}$}.
		\end{cases}
	\end{equation*}
	Here $E_{5}>0$ and $E_{6}>0$ are constants independent of $\mathcal{H}_1$, $\mathcal{H}_2$, and $n$.
	\item Let $Y$ be any complex Banach space and $U: \mathcal{B}(\mathcal{H}) \rightarrow Y$ be a non-null bounded linear operator with $\norm{U}<\lambda$. If $Y$ is a $s$-concave Banach lattice with $2 \leq s < \infty$ and there is a $1 \leq t <s$ such that $U$ is $(s,1)$-summing, then
	\begin{equation*}
		R_{\lambda}(B_Z,p,U) \geq D\, \left(\frac{\lambda^p - \norm{U}^p}{2\lambda^p - \norm{U}^p}\right)^{\frac{1}{p}} \, \left(\frac{\log \, n}{n}\right)^{1-\frac{1}{s}} 
	\end{equation*}
and \begin{equation*}
AP_{\lambda}(B_Z,p,U) \geq D\, \left(\frac{\lambda^p - \norm{U}^p}{2\lambda^p - \norm{U}^p}\right)^{\frac{1}{p}} \, \frac{\norm{Id:Z \rightarrow \ell^n _1}}{n}\,\left(\frac{\log \, n}{n}\right)^{1-\frac{1}{s}}
\end{equation*}
for some constant $D>0$ independent of $\mathcal{H}$ and $n$.
\end{enumerate}
	\end{thm}
\begin{rem}
\begin{enumerate}
	\item In particular, setting $q=1$ in Theorem \ref{thm-1.9} $(1)$ gives the following
		\begin{equation*}
		AP_{\lambda}(B_{Z}, p,U) \geq \frac{\Phi_{1}(n,\lambda,p,1)}{n \norm{Id:\ell^n _1 \rightarrow Z}}.
	\end{equation*}
	
	\item Since the lower bounds in $(1)$ of Theorem \ref{thm-1.9} are expressed in terms of the norms $\norm{Id:Z \rightarrow \ell^n _q}$ and $\norm{Id:\ell^n _q \rightarrow Z}$, varying $q$ over the interval $1 \leq q \leq \infty$ leads to an improvement of these estimates. The improved bounds are then given as follows:
	\begin{align*}
		R_{\lambda}(B_{Z}, p,U) \geq 
		\sup_{1\leq q \leq \infty} \left\{\frac{\Phi_{1}(n,\lambda,p,q)}{\norm{Id:Z \rightarrow \ell^n _q} \norm{Id:\ell^n _q \rightarrow Z}}\right\}
	\end{align*} 
	and
	\begin{equation*}
		AP_{\lambda}(B_{Z}, p,U) \geq 
		\sup_{1\leq q \leq \infty} \left\{\frac{\Phi_{1}(n,\lambda,p,q)\, \norm{Id:Z \rightarrow \ell^n _1}}{n\norm{Id:Z \rightarrow \ell^n _q} \norm{Id:\ell^n _q \rightarrow Z}}\right\}.
	\end{equation*}
\end{enumerate}
\end{rem}

Theorem \ref{thm-1.8-atm} provides a lower bound for the $n$-dimensional arithmetic Bohr radius $AP_\lambda(\Omega,p, U)$ in terms of the corresponding $n$-dimensional Bohr radius $R_\lambda(\Omega,p, U)$. For a given simply connected complete Reinhardt domain $\Omega \subset \mathbb{C}^n$, we now show that $AP_\lambda(\Omega,p, U)$ can be estimated in terms of the one-dimensional Bohr radius for $\mathbb{D}$, namely $R_\lambda(\mathbb{D},p,U)$.
\begin{thm}\label{tthm-1.10-atm}
Let $U:\mathcal{B}(\mathcal{H})\rightarrow Y$ be a non-null bounded liner operator and $\norm{U} < \lambda$.	Let $1\leq p<\infty$. Then for every $n$ we have 
	\begin{equation*}
		\frac{R_\lambda(\mathbb{D},p,U)}{n}\leq AP_\lambda(B_{Z},p, U)\leq \norm{Id: Z \rightarrow \ell^n_{1}}\, \frac{R_\lambda(\mathbb{D},p,U)}{n^{1/p}}.
	\end{equation*}
\end{thm}
By considering the function $\phi$ defined in \eqref{e-function} with $t=S(\Omega,B_{\ell^n_1})$ and proceeding along lines similar to those in the proof of Theorem \ref{tthm-1.10-atm}, one can verify that the result remains valid for any simply connected complete Reinhardt domain $\Omega$. We therefore state the result below and omit the details.
\begin{cor}\label{cor-1.10-atm}
	Let $U:\mathcal{B}(\mathcal{H})\rightarrow Y$ be a non-null bounded liner operator and $\norm{U} < \lambda$. Let $1\leq p<\infty$. Then for every $n$ we have 
	\begin{equation*}
		\frac{R_\lambda(\mathbb{D},p,U)}{n}\leq AP_\lambda(\Omega,p, U)\leq S(\Omega,B_{\ell^n_1})\, \frac{R_\lambda(\mathbb{D},p,U)}{n^{1/p}}.
	\end{equation*}
\end{cor}

In view of Theorem \ref{tthm-1.10-atm} and the fact $\norm{Id: B_{\ell^n_q}\rightarrow B_{\ell^n_1}} =n^{1-1/q}$ for $q \geq 1$, we deduce the following result for $Z=\ell^n_q$.
\begin{cor} \label{cor-1.5}
Let $U:\mathcal{B}(\mathcal{H})\rightarrow Y$ be a non-null bounded liner operator and $\norm{U} < \lambda$. Then for every $n$, $1\leq p<\infty$ and $1 \leq q \leq \infty$ we have 
\begin{equation*}
	\frac{R_\lambda(\mathbb{D},p,U)}{n}\leq AP_\lambda(B_{\ell^n_q},p, U)\leq  \frac{R_\lambda(\mathbb{D},p,U)}{n^{\frac{1}{p}+\frac{1}{q}-1}}.
\end{equation*}
\end{cor}
Unlike the Bohr radius constant $R_\lambda(\Omega,p, U)$, for which no exact value is currently known for any simply connected complete Reinhardt domain $\Omega \subset \mathbb{C}^n$ in dimension $n \geq 2$, the arithmetic Bohr radius can be explicitly computed in certain settings. In particular, for the unit ball of $\ell^n_1$ we are able to determine its precise value. This follows directly from Corollary \ref{cor-1.5} by taking $p=1$ and $q=1$.
\begin{cor} \label{cor-1.8}
Let $U:\mathcal{B}(\mathcal{H})\rightarrow Y$ be a non-null bounded liner operator and $\norm{U} < \lambda$.	Then for every $n\in \mathbb{N}$, we have 
	\begin{equation*}
		AP_\lambda(B_{\ell^n_1},1, U)=\frac{R_\lambda(\mathbb{D},1,U)}{n}.
	\end{equation*}
\end{cor}
It worth mentioning that Corollary \ref{cor-1.8} is actually pluriharmonic analogue of \cite[Theorem 3.1]{defant-2008}.
\vspace{1mm}

In view of the corollaries \ref{cor-1.10-atm}, \ref{cor-1.5}, and \ref{cor-1.8}, it is evident that to estimate $AP_\lambda(\Omega,p, U)$ it is enough to estimate $R_\lambda(\mathbb{D},p,U)$. To do this, we prove the following result in the case $U$ is the identity operator on $\mathbb{C}$ and $p=1$.
\begin{lem} \label{lem-1.1}
	
If $U=I:\mathbb{C} \rightarrow \mathbb{C}$ and $1 <\lambda$, then 
	\begin{equation*}
	R_\lambda(\mathbb{D},1,\mathbb{C}) \geq \frac{(\lambda-1)\pi}{4+(\lambda-1)\pi}.
	\end{equation*}
 \comment{If $U=I:\mathcal{B}(\mathcal{H}) \rightarrow \mathcal{B}(\mathcal{H})$ and $1 <\lambda$, then
 	\begin{equation*}
 		R_\lambda(\mathbb{D},1,\mathcal{B}(\mathcal{H})) \geq 	
 \end{equation*}}
\end{lem}
 Now, by making use of Lemma \ref{lem-1.1} and the corollaries \ref{cor-1.10-atm}, \ref{cor-1.5}, and \ref{cor-1.8}, we deduce the following important estimates, which we record as a corollary.
 \begin{cor}
 Let $\Omega \subset \mathbb{C}^n$ be a bounded simply connected complete Reinhardt domain. Let $U$ be the identity operator on $\mathbb{C}$ and $1 <\lambda$. Then for every $n \in \mathbb{N}$, we have
   
 	$$
 	\frac{(\lambda-1)\pi}{(4+(\lambda-1)\pi)\,n} \leq AP_\lambda(\Omega,1, \mathbb{C}).
 	$$
 \end{cor}
 This corollary actually provides a universal lower bound of the arithmetic Bohr radius for complex valued pluriharmonic functions on any bounded simply connected complete Reinhardt domain $\Omega \subset \mathbb{C}^n$.
 \vspace{1mm}
  
On the other hand, by applying the same arguments as in the pluriharmonic case in Theorem \ref{thm-1.8-atm}, Theorem \ref{tthm-1.10-atm}, and Corollary \ref{cor-1.10-atm}, we obtain the relation between the Bohr radius and the arithmetic Bohr radius for any Banach space $X$-valued holomorphic functions.
\begin{thm} \label{thm-1.8-atm-holo}
	Let $Z=(\mathbb{C}^n, ||.||)$ be a Banach space such that $\chi(\{e_k\}^n_{k=1})=1$. Then 
	\begin{equation*}
		A_\lambda(B_{Z},p,U) \geq \frac{\norm{Id: Z \rightarrow \ell^n_{1}}}{n} \, K_\lambda(B_{Z},p,U).
	\end{equation*}
\end{thm}

\begin{thm}\label{tthm-1.10-atm-holo}
	Let $U:X\rightarrow Y$ be a non-null bounded liner operator and $\norm{U} < \lambda$.	Let $1\leq p<\infty$. Then for every $n$ we have 
	\begin{equation*}
		\frac{K_\lambda(\mathbb{D},p,U)}{n}\leq A_\lambda(B_{Z},p, U)\leq \norm{Id: Z \rightarrow \ell^n_{1}}\, \frac{K_\lambda(\mathbb{D},p,U)}{n^{1/p}}.
	\end{equation*}
\end{thm}
Theorems \ref{thm-1.8-atm-holo} and \ref{tthm-1.10-atm-holo} remain valid for any bounded complete Reinhardt domain, not only for $B_Z$. Moreover, we recover \cite[Lemma 4.3]{defant-2008} from Theorem \ref{thm-1.8-atm-holo} by choosing $U=I:\mathbb{C} \rightarrow \mathbb{C}$ and $p=1$.
\par In particular, Theorem \ref{tthm-1.10-atm-holo} yields the following result.
\begin{cor} \label{cor-1.8-a}
	Let $U:X\rightarrow Y$ be a non-null bounded liner operator and $\norm{U} < \lambda$.	Then for every $n \in \mathbb{N}$ we have 
	\begin{equation*}
		A_\lambda(B_{\ell^n_1},1, U)=\frac{K_\lambda(\mathbb{D},1,U)}{n}.
	\end{equation*}
\end{cor}
It is interesting to note that \cite[Theorem 3.1]{defant-2008} follows from Corollary \ref{cor-1.8-a} by choosing $U$ the identity operator on $\mathbb{C}$.
\begin{rem}
	For $X=\mathcal{B}(\mathcal{H})$, the bounds for $AP_{\lambda}(B_{Z}, p,X)$ established in Theorems \ref{thm-1.2-atm}–\ref{thm-1.6} extend to the constant $A_{\lambda}(B_{Z}, p,X)$ for any complex Banach space $X$, possibly with different constants. In particular, for $X=\mathcal{B}(\mathcal{H})$, this implies that the asymptotic behaviors of $AP_{\lambda}(B_{Z}, p,X)$ and $A_{\lambda}(B_{Z}, p,X)$ coincide . Since the proofs follow the same lines as those of the corresponding results for $AP_{\lambda}(B_{Z}, p,X)$, we omit the details. Moreover, an analogue of Theorem \ref{tthm-1.8-atm} holds for the constant $A_{\lambda}(\Omega, p,X)$ for any bounded complete Reinhardt domain.
\end{rem}
	\section{Application}
	This section is devoted to examples that demonstrate the scope of our results when applied to various classes of sequence spaces. Our first example is the family of mixed Minkowski spaces, which form a natural generalization of the classical Minkowski spaces. Given $m,n \in \mathbb{N}$ and $1 \leq s,t \leq \infty$, the mixed Minkowski space is defined as 
	$$\ell^m_s(\ell^n_t):=\{(z_{k})^m_{k=1}:z_{1}, \ldots,z_{m} \in \mathbb{C}^n\}$$ 
	with $\norm{(z_{k})^m_{k=1}}_{s,t}:=(\sum_{k=1}^{m}\norm{z_{k}}^s_t)^{1/s}$. In the special case $s=t$, this construction coincides with the classical Minkowski space $\ell^{mn}_s$. While earlier contributions such as \cite{Allu-Pal-arithnmetic,A-H-P-banach} dealt primarily with classical Minkowski spaces and the holomorphic case, the framework proposed in this article extends beyond that setting. It treats a significantly larger family of sequence spaces and enables a systematic study of holomorphic as well as pluriharmonic functions.
	\begin{example}[{\bf Mixed spaces}]
		Let $m,n \geq 2$. Then
		\begin{enumerate}
			\item Let $1 \leq s ,t \leq 2$. If $X=\mathcal{B}(\mathcal{H})$ is finite dimensional, then
			\begin{equation*}
				AP_{\lambda}(B_{\ell^m_s(\ell^n_t)}, p,X) \geq	\begin{cases}
					E_{1}(X)\,  \, \max \left\{\frac{1}{n^{2/t}\,m^{2/s}}, \, \frac{1}{e\,mn }\right\}\, \left(\frac{\lambda^p -1}{2\lambda^p - 1}\right)^{\frac{1}{p}} & \text{for $p=1$}, \\[2mm]
					E_{3}\, n^{\frac{1}{p}-\frac{2}{tp}-\frac{1}{t}}\,m^{\frac{1}{p}-\frac{2}{sp}-\frac{1}{s}} \left(\frac{\lambda^p -1}{2\lambda^p - 1}\right)^{\frac{1}{p}}\, & \text{for $p \geq 2$}, \\[2mm]
					
					E_{2}(X)  \max \left\{\frac{1}{ n^{\frac{2}{t}+\frac{1}{p}-1}m^{\frac{2}{s}+\frac{1}{p}-1}},  \frac{m^{2+\frac{4}{ps}-\frac{3}{p}-\frac{4}{s}}}{e^{(1-\theta)} n^{\frac{3}{p}+\frac{4}{t} -2-\frac{4}{pt}}}\right\}\left(\frac{\lambda^p -1}{2\lambda^p - 1}\right)^{\frac{1}{p}} \hspace{-1mm} & \text{for $1<p< 2$}, 
				\end{cases}
			\end{equation*}
			and 
			$$AP_{\lambda}(B_{\ell^m_s(\ell^n_t)}, p,X) \leq d\,  \lambda^{\frac{1}{\log\,(mn)}}\,\, \sqrt{\log \, (mn)} \, \left(\frac{1}{mn}\right)^{\frac{1}{p}}.
			$$
			If $X=\mathcal{B}(\mathcal{H})$ is infinite dimensional then
			$$\frac{\tilde{\Psi}_{1}(mn,\lambda,p,2)}{m^{1/s -1/2}\,n^{1/t-1/2}} \, \leq AP_{\lambda}(B_{Z}, p,X) \leq  \tilde{\Psi}_{2}(mn,\lambda,p,2).
			$$ 
			\item Let $2 \leq s,t \leq \infty$. If $X=\mathcal{B}(\mathcal{H})$ is finite dimensional then
			\begin{equation*}
				AP_{\lambda}(B_{\ell^m_s(\ell^n_t)}, p,X) \geq	\begin{cases}
					E_{1}(X)\,  \, \max \left\{\frac{1}{n^{\frac{1}{2}+\frac{1}{t}}\, m^{\frac{1}{2}+\frac{1}{s}}}, \, \frac{1}{e\, mn}\right\}\, \left(\frac{\lambda^p -1}{2\lambda^p - 1}\right)^{\frac{1}{p}} & \text{for $p=1$}, \\[2mm]
					E_{3}\, n^{-\frac{1}{t}}\, m^{-\frac{1}{s}}\left(\frac{\lambda^p -1}{2\lambda^p - 1}\right)^{\frac{1}{p}}\, & \text{for $p \geq 2$}, \\[2mm]
					
					E_{2}(X)  \max \left\{\frac{1}{n^{\frac{1-\theta}{2}+\frac{1}{t}}\,m^{\frac{1-\theta}{2}+\frac{1}{s}}}, \frac{m^{\theta -\frac{\theta}{s}-1}}{e^{(1-\theta)} n^{1-\theta +\frac{\theta}{t}}}\right\}\left(\frac{\lambda^p -1}{2\lambda^p - 1}\right)^{\frac{1}{p}}  \hspace{-2mm} & \text{for $1<p< 2$}, 
				\end{cases}
			\end{equation*}
			and
			$$AP_{\lambda}(B_{\ell^m_s(\ell^n_t)}, p,X) \leq d\,n^{1-\frac{1}{t}-\frac{1}{p}}m^{1-\frac{1}{s}-\frac{1}{p}}\, \lambda^{\frac{2}{\log \, n}}\, \sqrt{\log mn}.
			$$
			If $X=\mathcal{B}(\mathcal{H})$ is infinite dimensional then
			$$\tilde{\Psi}_{1}(mn,\lambda,p,2)\, \leq AP_{\lambda}(B_{Z}, p,X) \leq m^{\frac{1}{2}-\frac{1}{s}}\,n^{\frac{1}{2}-\frac{1}{t}} \, \tilde{\Psi}_{2}(mn,\lambda,p,2).
			$$
		\end{enumerate}
		Here, the constants $E_{1}(X)$, $E_{2}(X)$, $E_{3}$, and $d$ are as in Theorem \ref{thm-1.2-atm}.	
	\end{example}
	\begin{pf}
	For $1\leq s,t,l,r \leq \infty$, it is easy to check that the identity operator between mixed sequence spaces factorizes as 
		$$\norm{Id:\ell^m_s(\ell^n_t) \rightarrow \ell^m_l(\ell^n_r)}= \norm{Id:\ell^m_s \rightarrow \ell^m_l}\, \norm{Id:\ell^n_t \rightarrow \ell^n_r}$$
		 (see \cite[p. 188]{defant-2003}). Moreover, for finite-dimensional $\ell^n_p$ spaces, the norm of the identity operator satisfies: $\norm{Id:\ell^n_s \rightarrow \ell^n_t}=1$ for $s \leq t$ and $\norm{Id:\ell^n_s \rightarrow \ell^n_t}=n^{1/t - 1/s}$ for $s>t$. Applying these observations together with Theorems \ref{thm-1.2-atm}, \ref{thm-1.4-atm}, and \ref{thm-1.6}, we obtain the desired result. This completes the proof.
	\end{pf}

	The scope of our results further includes classical symmetric Banach sequence spaces such as Lorentz spaces $\ell_{s,t}$ as well as Orlicz spaces $\ell_{\psi}$; see \cite{Lindenstrauss-book} for their definitions. 
	\begin{example}[{\bf Lorentz spaces}]
		Let $1 \leq s,t \leq \infty$.
		\begin{enumerate}
			\item Let $1 \leq s <2$ and $1\leq t \leq \infty$, or $s=2$ and $1 \leq t \leq 2$. If $X=\mathcal{B}(\mathcal{H})$ is finite dimensional 
			\begin{equation*}
				 AP_{\lambda}(\ell^n_{s,t}, p,X) \geq	\begin{cases}
					E_{1}(X)\,  \max \left\{\frac{1}{n^{\frac{1}{2}+\frac{1}{s}}\, \norm{Id: \ell^n_{2} \rightarrow \ell^n_{s,t}}}, \, \frac{1}{e\,n}\right\}\, \left(\frac{\lambda^p -1}{2\lambda^p - 1}\right)^{\frac{1}{p}} & \text{for $p=1$}, \\[2mm]
					E_{3}\, \frac{1}{n^{\frac{1}{s}}}\norm{Id:\ell^n _2 \rightarrow \ell^n_{s,t}}^{-\frac{2}{p}}\, \left(\frac{\lambda^p -1}{2\lambda^p - 1}\right)^{\frac{1}{p}}\, & \text{for $p \geq 2$}, \\[2mm]
					
					E_{2}(X)\,  \max \left\{\frac{1}{n^{1-\frac{\theta}{2}+\frac{1}{s}} \norm{Id: \ell^n_{2} \rightarrow \ell^n_{s,t}}},  \frac{\norm{Id:\ell^n _2 \rightarrow \ell^n_{s,t}}^{-\theta}}{e^{(1-\theta)}\,n^{1-\theta+\frac{\theta}{s}}}\right\}\left(\frac{\lambda^p -1}{2\lambda^p - 1}\right)^{\frac{1}{p}} \hspace{-2mm} & \text{for $1<p< 2$}, 
				\end{cases}
			\end{equation*}
			and 
			$$AP_{\lambda}(B_{\ell^n_{s,t}}, p,X) \leq d\,  \lambda^{\frac{1}{\log\,n}}\,\, \sqrt{\log \, n} \, \left(\frac{1}{n}\right)^{\frac{1}{p}}.
			$$
			If $X=\mathcal{B}(\mathcal{H})$ is infinite dimensional, then
			$$\frac{\tilde{\Psi}_{1}(n,\lambda,p,2)}{\norm{I:\ell^n _2 \rightarrow \ell^n_{s,t}}} \, \leq AP_{\lambda}(B_{\ell^n_{s,t}}, p,X)  \leq  \tilde{\Psi}_{2}(n,\lambda,p,2).
			$$
			\item Let $2<s \leq \infty$ and $1 \leq t \leq \infty$. If $X=\mathcal{B}(\mathcal{H})$ is finite dimensional
			\begin{equation*} 
			AP_{\lambda}(\ell^n_{s,t}, p,X) \geq	\begin{cases}
				E_{1}(X)\,  \max \left\{\frac{1}{n^{\frac{1}{2}+\frac{1}{s}}}, \, \frac{1}{e\,n}\right\}\, \left(\frac{\lambda^p -1}{2\lambda^p - 1}\right)^{\frac{1}{p}} & \text{for $p=1$}, \\[2mm]
				E_{3}\, \frac{1}{n^{\frac{1}{s}}}\, \left(\frac{\lambda^p -1}{2\lambda^p - 1}\right)^{\frac{1}{p}}\, & \text{for $p \geq 2$}, \\[2mm]
				
				E_{2}(X)\,  \max \left\{\frac{1}{n^{1-\frac{\theta}{2}+\frac{1}{s}}},  \frac{1}{e^{(1-\theta)}\,n^{1-\theta+\frac{\theta}{s}}}\right\}\left(\frac{\lambda^p -1}{2\lambda^p - 1}\right)^{\frac{1}{p}} \hspace{-2mm} & \text{for $1<p< 2$}, 
			\end{cases}
			\end{equation*}
			and $$AP_{\lambda}(B_{\ell^n_{s,t}}, p,X) \leq d\,  \lambda^{\frac{1}{\log\,n}}\,\, \sqrt{\log \, n} \, \left(\frac{1}{n}\right)^{\frac{1}{p}+\frac{1}{s}-\frac{1}{2}}.
			$$ If $X=\mathcal{B}(\mathcal{H})$ is infinite dimensional, then $$\tilde{\Psi}_{1}(n,\lambda,p,2) \, \leq AP_{\lambda}(B_{\ell^n_{s,t}}, p,X)  \leq \left(\frac{1}{n}\right)^{\frac{1}{s}-\frac{1}{2}} \, \tilde{\Psi}_{2}(n,\lambda,p,2).
			$$
		\end{enumerate}
	\end{example}
	\begin{pf}
		We begin by observing that $\norm{\sum_{k=1}^{n}e^{*}_{k}}_{\ell_{s,t}^{*}}=n^{1-1/s}$. Moreover, Lorentz sequence spaces satisfy the following lexicographical inclusion property: $$\ell_{p,q} \subseteq \ell_{s,t}\,\, \mbox{if and only if}\,\,\, p<s, \,\, \mbox{or}\,\, p=s\,\, \mbox{and}\,\, q \leq t.
		$$
		 These facts, together with Theorems \ref{thm-1.3-a}(1), \ref{thm-1.4-atm}, and \ref{thm-1.6}, immediately yield the assertion in part $(1)$. 
		 \vspace{1mm}
		 
		 We now turn to part $(2)$. Suppose first that $2<s \leq \infty$ and $2 \leq t \leq \infty$. In this case, $\ell_{s,t}$ is $2$-convex (see \cite[p.189]{defant-2003}), and the conclusion follows directly from Theorems \ref{thm-1.3-a}$(2)$, \ref{thm-1.4-atm}, and \ref{thm-1.6}. It remains to consider the case $2 < s \leq \infty$, $1 \leq t \leq 2$. Here we apply Theorems \ref{thm-1.2-atm}, \ref{thm-1.4-atm}, and \ref{thm-1.6}, together with the estimates 
		 $$
		 \norm{Id: \ell^n_{2} \rightarrow \ell^n_{s,t}} \leq 1\,\,\,\mbox{and}\,\,\, \norm{Id: \ell^n_{s,t} \rightarrow \ell^n_{2}} \leq n^{1/2 - 1/s}
		 $$ 
		 (see \cite[p.189]{defant-2003}). This completes the proof.
	\end{pf}
	\begin{example} [{\bf Orlicz spaces}]
		Let $\psi$ be an Orlicz function which satisfies the $\Delta_{2}$-condition. 
		\begin{enumerate}
			\item Let $a^2 \leq T \psi(a)$ for all $a$ and some $T>0$.  If $X=\mathcal{B}(\mathcal{H})$ is finite dimensional, then
			\begin{equation*}
				AP_{\lambda}(B_{\ell^n_{\psi}}, p,X) \geq	\begin{cases}
					E_{1}(X)\,  \, \max \left\{\frac{\psi^{-1}(1/n)}{\sqrt{n}\, \norm{Id: \ell^n_{2} \rightarrow \ell^n_{\psi}}}, \, \frac{1}{e\,n}\right\}\, \left(\frac{\lambda^p -1}{2\lambda^p - 1}\right)^{\frac{1}{p}} & \text{for $p=1$}, \\[2mm]
					E_{3}\, \psi^{-1}(1/n)\,\norm{Id:\ell^n _2 \rightarrow \ell^n_{\psi}}^{-\frac{2}{p}}\, \left(\frac{\lambda^p -1}{2\lambda^p - 1}\right)^{\frac{1}{p}}\, & \text{for $p \geq 2$}, \\[2mm]
					
					E_{2}(X)  \max \left\{\frac{\psi^{-1}(1/n)}{(\sqrt{n})^{1-\theta} \norm{Id: \ell^n_{2} \rightarrow \ell^n_{\psi}}},  \frac{\norm{Id:\ell^n _2 \rightarrow \ell^n_{\psi}}^{-\theta}}{(en)^{(1-\theta)} (\psi^{-1}(1/n))^{-\theta}}\right\}\left(\frac{\lambda^p -1}{2\lambda^p - 1}\right)^{\frac{1}{p}} \hspace{-4mm}  & \text{for $1<p< 2$}, 
				\end{cases}
			\end{equation*}
			and $$AP_{\lambda}(B_{\ell^n_{\psi}}, p,X) \leq d\,  \lambda^{\frac{1}{\log\,n}}\,\, \sqrt{\log \, n} \, \left(\frac{1}{n}\right)^{\frac{1}{p}}.
			$$
			If $X=\mathcal{B}(\mathcal{H})$ is infinite dimensional, 
			$$\frac{\tilde{\Psi}_{1}(n,\lambda,p,2)}{\norm{Id:\ell^n _2 \rightarrow \ell^n_{\psi}}} \, \leq AP_{\lambda}(B_{\ell^n_{\psi}}, p,X) \leq  \tilde{\Psi}_{2}(n,\lambda,p,2).
			$$
			
			\item Let $\psi(\beta \,a) \leq T \beta^2\,\psi(a)$ for $0 \leq \beta,a \leq 1$ and some $T>0$. If $X=\mathcal{B}(\mathcal{H})$ is finite dimensional,
			\begin{equation*}
				AP_{\lambda}(B_{\ell^n_{\psi}}, p,X) \geq	\begin{cases}
					E_{1}(X)\,  \, \max \left\{\frac{\psi^{-1}(1/n)}{\sqrt{n}}, \, \frac{1}{e\,n}\right\}\, \left(\frac{\lambda^p -1}{2\lambda^p - 1}\right)^{\frac{1}{p}} & \text{for $p=1$}, \\[2mm]
					E_{3}\,\psi^{-1}(1/n)\, \left(\frac{\lambda^p -1}{2\lambda^p - 1}\right)^{\frac{1}{p}}\, & \text{for $p \geq 2$}, \\[2mm]
					
					E_{2}(X)\,  \max \left\{\frac{\psi^{-1}(1/n)}{(\sqrt{n})^{1-\theta}}, \, \frac{1}{(en)^{(1-\theta)}\, (\psi^{-1}(1/n))^{-\theta}}\right\}\left(\frac{\lambda^p -1}{2\lambda^p - 1}\right)^{\frac{1}{p}}  & \text{for $1<p< 2$}, 
				\end{cases}
			\end{equation*}
			and 
			$$AP_{\lambda}(B_{\ell^n_{\psi}}, p,X) \leq d\,  \lambda^{\frac{1}{\log\,n}}\, \,\psi^{-1}\left(\frac{1}{n}\right)\,\, \sqrt{\log \, n}\, \left(\frac{1}{n}\right)^{\frac{1}{p}-\frac{1}{2}}.$$
			If $X=\mathcal{B}(\mathcal{H})$ is infinite dimensional, 
			$$\tilde{\Psi}_{1}(n,\lambda,p,2) \, \leq AP_{\lambda}(B_{\ell^n_{\psi}}, p,X) \leq \sqrt{n}\, \psi^{-1}\left(\frac{1}{n}\right) \, \tilde{\Psi}_{2}(n,\lambda,p,2).
			$$
		\end{enumerate}
	\end{example}
	\begin{pf}
		It is important to note the facts $\norm{\sum_{k=1}^{n}e^{*}_{k}}_{\ell_{\psi}^{*}}\, \norm{\sum_{k=1}^{n}e_{k}}_{\ell_{\psi}}=n$ 
		and
		$$
		\norm{\sum_{k=1}^{n}e_{k}}_{\ell_{\psi}}=\frac{1}{\psi^{-1}\left(\frac{1}{n}\right)}
		$$
		(see \cite[p.192]{defant-2003}). The condition imposed on $\psi$ in part $(1)$ guarantees that $\ell_{\psi} \subseteq \ell_2$.  On the other hand, the condition in part $(2)$ ensures that $\ell_{\psi}$ is $2$-convex (see \cite[p.189]{defant-2003}). Taking these properties into account, the desired conclusion follows directly from Theorems \ref{thm-1.3-a} and \ref{thm-1.6}, together with the preceding facts. This completes the proof.
	\end{pf}	
	\section{Preliminaries}
	We adopt standard notation and terminology from Banach space theory. Throughout, all Banach spaces $W$ are assumed to be complex. Their topological duals are denoted by $W^{*}$, and their open unit balls by $B_{W}$. A Banach space $W$ is said to have cotype $t\in [0,\infty]$ if there exists a constant $C>0$ such that, for every finite family vectors $w_{1},\ldots,w_{n}\in W$, 
	$$\left(\sum_{k=1}^{n}\norm{w_{k}}^t\right)^{\frac{1}{t}} \leq C \left(\int_{0}^{1}\norm{\sum_{k=1}^{n}r_{k}(s)w_{k}}^2\, ds\right)^{\frac{1}{2}},$$
	where $r_{k}$ denotes the $k$-th Rademacher function on $[0,1]$. We define 
	$$Cot(W):=\inf\{2 \leq t \leq \infty: W\, \mbox{has cotype}\, t\}.
	$$
	It is well known that every Banach space has cotype $\infty$.  When $Cot(W)=\infty$,  we adopt the convention $1/Cot(W)=0$. A Schauder basis $\{w_{k}\}$ of a Banach space $W$ is called unconditional if there is a constant $c\geq 1$ such that $$\norm{\sum_{j=1}^{k}\varepsilon_{j}\mu_{j}w_{j}} \leq c \norm{\sum_{j=1}^{k}\mu_{j}w_{j}}
	$$
	for all $k\in\mathbb{N}$, all scalars $\mu_{j}\in\mathbb{C}$, and all $\varepsilon_{j}\in\mathbb{C}$ with $|\varepsilon_{j}|\le 1$. The smallest such constant is denoted by $\chi({w_{k}})$ and is called the unconditional basis constant of ${w_{k}}$. We say that ${w_{k}}$ is a $1$-unconditional basis if $\chi({w_{k}})=1$. The unconditional basis constant of the Banach space $W$ is defined by $\chi(W):=\inf \chi(\{w_{k}\})\in [1,\infty]$, the infimum taken over all unconditional basis $\{w_{k}\}$ of $W$. A Banach space $Z$ satisfying $\ell_1 \subset Z \subset c_{0}$ (with normal inclusions) is called a Banach sequence space if the canonical unit vectors ${e_{k}}$ form a $1$-unconditional basis of $Z$. A Banach lattice $Z$ is said to be $2$-convex if there exists a constant $C>0$ such that 
	$$\norm{\left(\sum_{k=1}^{n}|x_{k}|^2\right)^{\frac{1}{2}}} \leq C 	\left(\sum_{k=1}^{n}\norm{x_{k}}^2\right)^{\frac{1}{2}}
	$$
	for all finite families $x_{1},\ldots,x_{n}\in Z$. Let $Z=(\mathbb{C}^n, ||.||)$ be a Banach space, $Y$ any Banach space, and $m \in \mathbb{N}$. We denote by $\mathcal{P}(^m Z,Y)$ the space of all $m$-homogeneous polynomials $Q:Z \rightarrow Y$ of the form $Q(z)=\sum_{|\alpha|=m}c_{\alpha}z^{\alpha}$, together with the norm $\norm{Q}_{\mathcal{P}(^m Z,Y)}:=\sup_{z \in B_{Z}}||Q(z)||_{Y}$. The unconditional basis constant of the family of monomials $z^{\alpha}$, $\alpha \in (\mathbb{N} \cup \{0\})^n$ is denoted by $\chi_{M}(\mathcal{P}(^m Z,Y))$.
 We denote by $\mathcal{PH}(^m Z,\mathcal{B}(H))$ the space of all $m$-homogeneous pluriharmonic polynomials $P:Z \rightarrow \mathcal{B}(H)$ of the form 
	\begin{equation} \label{e-1.3-aa}
		P(z)= \sum_{|\alpha|=m} a_{\alpha}\, z^{\alpha} +  \sum_{|\alpha|=m} b^{*}_{\alpha}\, \bar{z}^{\alpha}
	\end{equation} 
	and set $\norm{P}_{\mathcal{PH}(^m Z,\mathcal{B}(H))}:=\sup_{z \in B_{Z}}|||P(z)||_{\mathcal{B}(H)}$. 
	\par A bounded linear operator $U:W\rightarrow Y$ between two Banach spaces $W$ and $Y$ is said to be $(r,t)$-summing, $r,t \in [1,\infty)$, if exists a constant $C>0$ such that for every finite collection $w_{1}, \ldots,w_{n} \in W$,
	$$\left(\sum_{k=1}^{n} \norm{U(w_{k})}^r\right)^{1/r} \leq C \, \sup_{\phi \in B_{X^{*}}} \left(\sum_{k=1}^{n} |\phi(w_{k})|^t\right)^{1/t}.
	$$
	The smallest constant $C$ for which this inequality holds is denoted by $\pi_{r,t}(U)$. In the special case $r=t$, the operator $U$ is called $r$- summing, and we write $\pi_{r}(U)$ instead of $\pi_{r,r}(U)$.
	\section{Auxiliary Lemmas} 
	The following definition introduces the $m$-homogeneous analogue of the definition of Bohr radius for pluriharmonic functions, which plays a fundamental role in the analysis carried out in this article. Let $\Omega\subset \mathbb{C}^n$ be a simply connected complete Reinhardt domain and $n\in \mathbb{N}$. Let $U:X\rightarrow Y$ be a bounded liner operator and $\norm{U} \leq \lambda$. For $1 \leq p < \infty$, the $\lambda$-powered Bohr radius of $U$ with respect to $\mathcal{PH}(^m \Omega,X)$, denoted by $R^m_{\lambda}(\Omega, p,U)$, is defined to be the supremum of all $r\geq 0$ such that for every pluriharmonic polynomials $P:\Omega \rightarrow X$ of the form \eqref{e-1.3-aa} we have 
	\begin{equation} \label{e-1.4-aa}
		\sup_{z \in r\Omega}\, \sum_{|\alpha|=m} (\norm{U(a_{\alpha})}^p_{Y} + \norm{U(b_{\alpha})}^p_{Y})|z^\alpha|^p \leq \lambda^p\,\norm{f}^p_{\Omega,X}. 
	\end{equation}
	
	\noindent Here $\norm{f}_{\Omega,X}:=\sup_{z \in \Omega}\,\norm{f(z)}_{X}$. Set $R^m(\Omega, p,U):=R^m_{1}(\Omega, p,U)$, $R^m_{\lambda}(\Omega, p,X):=R^m_{\lambda}(\Omega, p,U)$ whenever $U=I:X\rightarrow X$, $R^m(\Omega, p,X):=R^m_{1}(\Omega, p,X)$, $R^m_{\lambda}(\Omega, p):=R^m_{\lambda}(\Omega, p,\mathbb{C})$, and $R^m(\Omega,p):=R^m_1(\Omega,p)$. It is straightforward that $R^m_{\lambda}(\Omega, p,U)=\lambda^{1/m}\,R^m(\Omega, p,U)$. Similarly as above, we can define the constant $K^m_{\lambda}(\Omega, p,U)$ for the class $\mathcal{P}^n_m(W)$ {\it i.e.}, for any Banach space $W$ valued $m$-homogeneous holomorphic polynomial $Q(z)=\sum_{|\alpha|=m}c_{\alpha}\, z^{\alpha}$. Clearly, we have $K^m_{\lambda}(\Omega, p,U)=\lambda^{1/m}\,K^m(\Omega, p,U)$. Using Bohnenblust-Hille inequality, Defant {\it et al.} \cite{defant-2011} have shown that the constant $K^m(\mathbb{D}^n,1)$ is hypercontractive.
	\vspace{2mm}
	
The following lemma, proved in \cite{Himadri-local-Banach-1}, provides a coefficient-type Schwarz–Pick inequality for operator valued pluriharmonic functions defined on general complete Reinhardt domains. It plays an essential role in the proofs of several of our main results.
	\begin{lem} \label{lem-3.1} \cite{Himadri-local-Banach-1}
		Let $f \in \mathcal{PH}(Z,X)$ be of the form \eqref{e-1.3-a}. Then for all $|\alpha|=m\geq 1$, we have $\norm{\sum_{|\alpha|=m}(a_{\alpha} \pm b_{\alpha})\, z^{\alpha}}_{B_{Z},X} \leq 4 \,\norm{\norm{f}_{B_{Z},X}\,I-\real(a_{0})}$. Moreover, if $B_{Z}=B_{\ell^n_{q}}$ then $\norm{a_{\alpha}+b_{\alpha}}	\leq \frac{4}{\pi} \, \rho_{\alpha}\, \norm{\sum_{|\alpha|=m}(a_{\alpha}+b_{\alpha})\, z^{\alpha}}_{B_{Z},X}$,  $\norm{a_{\alpha}-b_{\alpha}}	\leq \frac{4}{\pi}\,\rho_{\alpha}\, \norm{\sum_{|\alpha|=m}(a_{\alpha}-b_{\alpha})\, z^{\alpha}}_{B_{Z},X}$, where $\rho_{\alpha}:=\left(\frac{|\alpha|^{|\alpha|}}{\alpha^{\alpha}}\right)^{1/q}$.
		\comment{\begin{enumerate}
				\item $\norm{\sum_{|\alpha|=m}(a_{\alpha}+b_{\alpha})\, z^{\alpha}}_{B_{Z},X} \leq 4 \,\gamma_{0} \norm{f}_{B_{Z},X}$,
				$\norm{\sum_{|\alpha|=m}(a_{\alpha}-b_{\alpha})\, z^{\alpha}}_{B_{Z},X} \leq 4 \, \gamma_{0}\, \norm{f}_{B_{Z},X}$;
				\item $\norm{a_{\alpha}+b_{\alpha}}	\leq \frac{4}{\pi} \, \rho_{\alpha}\, \norm{\sum_{|\alpha|=m}(a_{\alpha}+b_{\alpha})\, z^{\alpha}}_{B_{Z},X}$,  $\norm{a_{\alpha}-b_{\alpha}}	\leq \frac{4}{\pi}\,\rho_{\alpha}\, \norm{\sum_{|\alpha|=m}(a_{\alpha}-b_{\alpha})\, z^{\alpha}}_{B_{Z},X}$,
			\end{enumerate}
			where $\gamma_{0}:=\norm{I-\real(a_{0})}$ and $\rho_{\alpha}:=\left(\frac{|\alpha|^{|\alpha|}}{\alpha^{\alpha}}\right)^{\frac{1}{q}}$.}
	\end{lem}

	The following lemma, established in \cite{Himadri-local-Banach-1}, provides a fundamental connection between the Bohr radius for pluriharmonic functions and that for homogeneous pluriharmonic functions. Conceptually, it may be viewed as a pluriharmonic counterpart of \cite[Lemma 3.2]{defant-2012}, where analogous ideas were developed for Minkowski spaces in the setting of the powered Bohr inequality. In contrast to that earlier result, the following lemma is formulated in the general framework of Banach spaces.
	\begin{lem} \cite{Himadri-local-Banach-1} \label{lem-3.3}
		Let $X=\mathcal{B}(\mathcal{H})$ and $Y$ any complex Banach space, and $U: X \rightarrow Y$ be a non-null bounded linear operator such that $\norm{U}< \lambda$. Then, for all $p\in[1,\infty)$, $\lambda>1$, and $n \in \mathbb{N}$, we have
		\begin{enumerate}
			\item $\left(\frac{\lambda^p - \norm{U}^p}{2\lambda^p - \norm{U}^p}\right)^{\frac{1}{p}}\, \inf_{m\in \mathbb{N}}\, R^{m}_{\lambda}(B_{Z}, p,U) \leq R_{\lambda}(B_{Z}, p,U) \leq \inf_{m\in \mathbb{N}}\, R^{m}_{\lambda}(B_{Z}, p,U)$
			
			\item $\left(\frac{\lambda^p - \norm{U}^p}{\lambda^p - \norm{U}^p +1}\right)^{\frac{1}{p}}\, \inf_{m\in \mathbb{N}}\, R^{m}(B_{Z}, p,U) \leq R_{\lambda}(B_{Z}, p,U) \leq \lambda \, \inf_{m\in \mathbb{N}}\, R^{m}(B_{Z}, p,U)$.
		\end{enumerate}
	\end{lem}
	It worth mentioning that analogous results of Lemma \ref{lem-3.3} hold for the constants $K^m_{\lambda}(B_{Z}, p,U)$ and $K_{\lambda}(B_{Z}, p,U)$.
	\vspace{1mm}
	
The following lemma compares the arithmetic Bohr radius constants for two complete Reinhardt domains. This result serves as a key ingredient in the proofs of several of our main theorems.
\begin{lem} \label{lem-3.6}
	Let $\Omega_{1}$ and $\Omega_{2}$ be two bounded simply connected complete Reinhardt domains in $\mathbb{C}^n$. Then we have $$\frac{AP_{\lambda}(\Omega_{2}, p,U)}{S(\Omega_{2},\Omega_{1})} \, \leq AP_{\lambda}(\Omega_{1}, p,U) \leq S(\Omega_{1},\Omega_{2}) \, AP_{\lambda}(\Omega_{2}, p,U).
	$$
\end{lem}
\begin{pf}
Let $X=\mathcal{B}(\mathcal{H})$. We shall first prove the left hand inequality. Let $r\in \mathbb{R}^n_{\geq 0}$ be such that for all $f \in \mathcal{PH}(\Omega_2,X)$ of the form \eqref{e-1.3-a},
\begin{equation} \label{ee-1.8-atm-added}
	\sum_{m=0}^{\infty}	\sum_{|\alpha|=m} (\norm{U(a_{\alpha})}^p_{Y} + \norm{U(b_{\alpha})}^p_{Y})r^{p \alpha} \leq \lambda^{p} \norm{f}^p_{\Omega_2, X}.	
\end{equation}
Then, for every $\epsilon >0$, we obtain
\begin{align*}
&\sum_{m=0}^{\infty}	\sum_{|\alpha|=m} (\norm{U(a_{\alpha})}^p_{Y} + \norm{U(b_{\alpha})}^p_{Y}) \left(\frac{r}{S(\Omega_{2},\Omega_{1}) +\epsilon}\right)^{p \alpha} \\
& \leq \lambda^p \norm{\sum_{m=0}^{\infty} \sum_{|\alpha|=m} a_{\alpha}\, \left(\frac{z}{S(\Omega_{2},\Omega_{1}) +\epsilon}\right)^{\alpha} + \sum_{m=1}^{\infty} \sum_{|\alpha|=m} b^{*}_{\alpha}\, \bar{\left(\frac{z}{S(\Omega_{2},\Omega_{1}) +\epsilon}\right)}^{\alpha}}_{\Omega_2,X} \\
& \leq \norm{f}_{\Omega_{1},X}.
\end{align*} 
This implies that
\begin{equation*}
\frac{1}{S(\Omega_{2},\Omega_{1}) +\epsilon} \, \frac{\sum_{i=1}^{n} r_i}{n}=\frac{1}{n}\left(\frac{r_1}{S(\Omega_{2},\Omega_{1}) +\epsilon} +\cdots + \frac{r_n}{S(\Omega_{2},\Omega_{1}) +\epsilon}\right) \leq AP_{\lambda}(\Omega_{1}, p,U)
\end{equation*}
for all $\epsilon>0$.
Consequently, 
\begin{equation*}
\frac{1}{S(\Omega_{2},\Omega_{1})}.\, AP_{\lambda}(\Omega_{2}, p,U) \leq  AP_{\lambda}(\Omega_{1}, p,U). 
\end{equation*}
By interchanging $\Omega_{1}$ and $\Omega_{2}$ in the above inequality, we establish the right-hand inequality.
\end{pf}

	We recall the following remarkable result by Maurey and Pisier to prove the desired upper bound in the case of infintedimensional complex Banach space. 
	\begin{customthm}{A} \cite[Theorem 14.5]{diestel-abs-summing-1995} \label{thm-A}
		Given any infinite dimensional complex Banach space $X$, there exist $x_{1}, \ldots , x_{n} \in X$ for each $n \in \mathbb{N}$ such that $\norm{z}_{\infty}/2 \leq \norm{\sum_{j=1}^{n} x_{j}z_{j}} \leq \norm{z}_{Cot(X)}$  for every choice of $z=(z_{1}, \ldots,z_{n}) \in \mathbb{C}^n$. Clearly, setting $z=e_{j}$, gives $\norm{x_{j}}\geq 1/2$, where $e_{j}$ is the $j$-th canonical basis vector of $\mathbb{C}^n$.
	\end{customthm}
	\section{Proof of the main results}
	\begin{pf} [{\bf Proof of Theorem \ref{thm-1.8-atm}}]
		Since $\norm{Id: Z \rightarrow \ell^n_{1}}=\sup_{z \in B_{Z}} \norm{z}_{1}$, for any $0< \epsilon < R_\lambda(B_{Z},p,U)$, we can find an element $\tilde{z}\in B_{Z}$ such that 
		\begin{equation*}
			\norm{\tilde{z}}_{\ell^n_1}\geq \norm{Id: Z \rightarrow \ell^n_{1}}-\epsilon.
		\end{equation*}
		Let $s:=R_\lambda(B_{Z},p,U)-\epsilon$, $w:=s\tilde{z}$, and $r:=s|\tilde{z}|=|w|.$ Since $w\in s\, B_{Z}$ and $s<R_\lambda(B_{Z},p,U)$, for $f \in \mathcal{PH}(B_{z},\mathcal{B}(\mathcal{H}))$ be of the form \eqref{e-1.3-a}, we have
		\begin{align*}
			\sum_{m=0}^{\infty} \sum_{|\alpha|=m} (\norm{U(a_{\alpha})}^p_{Y} + \norm{U(b_{\alpha})}^p_{Y})r^{p \alpha} &=\sum_{m=0}^{\infty} \sum_{|\alpha|=m} (\norm{U(a_{\alpha})}^p_{Y} + \norm{U(b_{\alpha})}^p_{Y})|w^\alpha|^p \\
			&\leq \sup_{z \in s\, B_{Z}}\sum_{m=0}^{\infty} \sum_{|\alpha|=m} (\norm{U(a_{\alpha})}^p_{Y} + \norm{U(b_{\alpha})}^p_{Y})|z^\alpha|^p \\
			&\leq \lambda^{p} \norm{f}^p_{B_{Z}, \mathcal{B}(\mathcal{H})}.
		\end{align*}
		Therefore, we obtain 
		\begin{equation*}
			AP_\lambda(B_{Z},p,U) \geq \frac{1}{n}\sum_{i=1}^{n}r_i = \frac{R_\lambda(B_{Z},p,U)-\epsilon}{n}\norm{\tilde{z}}_{\ell^n_{1}}\geq \frac{R_\lambda(B_{Z},p,U)-\epsilon}{n}\left(\norm{Id: Z \rightarrow \ell^n_{1}}-\epsilon\right)
		\end{equation*}
		holds for all $\epsilon>0$, and so \eqref{e-1.1-added} immediately follows. 
		On the other hand, when $\norm{U}< \lambda$, the following lower bound of $R_\lambda(B_{Z},p,U)$ was established in \cite[Theorem 1.1]{Himadri-local-Banach-1}:
		\begin{equation*}
			R_{\lambda}(B_{Z}, p,U) \geq D. \frac{1}{\sup_{z \in B_{Z}}\norm{z}_{p}},
		\end{equation*}
	where $D$ as in the statement of Theorem \ref{thm-1.8-atm}.
		 Consequently, the second part of the theorem follows from \eqref{e-1.1-added} together with the preceding bound. This completes the proof.	
	\end{pf}
\begin{pf} [{\bf Proof of Theorem \ref{tthm-1.8-atm}}]
	We shall first prove the left-hand inequality in $(1)$. Let $r=(r_1,.\,.\,., r_n)\in \mathbb{R}^n_{\geq 0}$ be such that 
	\begin{equation*}
		\sum_{|\alpha|=m} (\norm{U(c_{\alpha})}^p_{Y} + \norm{U(d_{\alpha})}^p_{Y})r^{p \alpha} \leq \lambda^{p} \norm{g}^p_{\Omega, X}	
	\end{equation*}
	for all $g(z)=\sum_{|\alpha|=m} c_{\alpha}\, z^{\alpha} +  \sum_{|\alpha|=m} d^{*}_{\alpha}\, \bar{z}^{\alpha}\in \mathcal{PH}(^m\Omega,X)$. Fix $f\in \mathcal{PH}(\Omega,X)$ be of the form \eqref{e-1.3-a}. Let $R \in (0,1)$. Then 
	\begin{align} \label{ee-1.5-atm}
		&\sum_{m=0}^{\infty} \sum_{|\alpha|=m} (\norm{U(a_{\alpha})}^p_{Y} + \norm{U(b_{\alpha})}^p_{Y})(Rr)^{p \alpha}  \nonumber\\
		&= \norm{U(a_{0})}^p +  \sum_{m=1}^{\infty} R^{pm}\sum_{|\alpha|=m} (\norm{U(a_{\alpha})}^p_{Y} + \norm{U(b_{\alpha})}^p_{Y})r^{p \alpha} \nonumber\\ 
		& \leq \norm{U}^p \norm{a_{0}}^p +  \sum_{m=1}^{\infty} R^{pm} \lambda^{p} \norm{\sum_{|\alpha|=m} a_{\alpha}\, z^{\alpha} +  \sum_{|\alpha|=m} b^{*}_{\alpha}\, \bar{z}^{\alpha}}^p_{\Omega,X} \\
		&\leq \norm{U}^p \norm{f}^p _{\Omega,X} + \sum_{m=1}^{\infty} R^{pm} \lambda^{p} \norm{f}^p _{\Omega,X}= \left(\norm{U}^p + \lambda^p\, \frac{R^p}{1-R^p}\right) \norm{f}^p _{\Omega,X} \leq \lambda^p \norm{f}^p _{\Omega,X} ,\nonumber
	\end{align}
	provided 
	$$R \leq \left(\frac{\lambda^p - \norm{U}^p}{2\lambda^p - \norm{U}^p}\right)^{\frac{1}{p}}.$$ This yields
	\begin{equation*}
		\left(\frac{\lambda^p - \norm{U}^p}{2\lambda^p - \norm{U}^p}\right)^{\frac{1}{p}} \, \frac{1}{n} \, \sum_{i=1}^{n} r_{i} \leq AP_{\lambda}(B_{Z}, p,U),
	\end{equation*}
	and hence, the desired bound follows. The uper bound of the right-hand inequality of $(1)$ directly follows from the fact $\cup^{\infty}_{m=1} \mathcal{PH}(^m \Omega,X) \subseteq \mathcal{PH}(\Omega,X)$.
	\vspace{1mm}
	
	 We next prove the left-hand inequality in $(2)$. Let $r=(r_1,.\,.\,., r_n)\in \mathbb{R}^n_{\geq 0}$ be such that 
	\begin{equation*}
		\sum_{|\alpha|=m} (\norm{U(c_{\alpha})}^p_{Y} + \norm{U(d_{\alpha})}^p_{Y})r^{p \alpha} \leq \norm{g}^p_{\Omega, X}	
	\end{equation*}
	for all $g(z)=\sum_{|\alpha|=m} c_{\alpha}\, z^{\alpha} +  \sum_{|\alpha|=m} d^{*}_{\alpha}\, \bar{z}^{\alpha}\in \mathcal{PH}(^m\Omega,X)$. Fix $f\in \mathcal{PH}(\Omega,X)$ be of the form \eqref{e-1.3-a}. Let $R \in (0,1)$. Then 
	\begin{align} \label{ee-1.6-atm}
		&\sum_{m=0}^{\infty} \sum_{|\alpha|=m} (\norm{U(a_{\alpha})}^p_{Y} + \norm{U(b_{\alpha})}^p_{Y})(Rr)^{p \alpha}  \nonumber\\
		&= \norm{U(a_{0})}^p +  \sum_{m=1}^{\infty} R^{pm}\sum_{|\alpha|=m} (\norm{U(a_{\alpha})}^p_{Y} + \norm{U(b_{\alpha})}^p_{Y})r^{p \alpha} \nonumber\\ 
		& \leq \norm{U}^p \norm{a_{0}}^p +  \sum_{m=1}^{\infty} R^{pm}  \norm{\sum_{|\alpha|=m} a_{\alpha}\, z^{\alpha} +  \sum_{|\alpha|=m} b^{*}_{\alpha}\, \bar{z}^{\alpha}}^p_{\Omega,X} \\
		&\leq \norm{U}^p \norm{f}^p _{\Omega,X} + \sum_{m=1}^{\infty} R^{pm} \norm{f}^p _{\Omega,X}= \left(\norm{U}^p + \frac{R^p}{1-R^p}\right) \norm{f}^p _{\Omega,X} \leq  \lambda^p \norm{f}^p _{\Omega,X} ,\nonumber
	\end{align}
	provided 
	$$R \leq \left(\frac{\lambda^p - \norm{U}^p}{\lambda^p - \norm{U}^p +1}\right)^{\frac{1}{p}}.$$ This yields
	\begin{equation*}
		\left(\frac{\lambda^p - \norm{U}^p}{\lambda^p - \norm{U}^p+1}\right)^{\frac{1}{p}} \, \frac{1}{n} \, \sum_{i=1}^{n} r_{i} \leq AP_{\lambda}(B_{Z}, p,U),
	\end{equation*}
	and hence, the desired bound follows. The uper bound of the right-hand inequality in $(2)$ directly follows from the fact $\cup^{\infty}_{m=1} \mathcal{PH}(^m \Omega,X) \subseteq \mathcal{PH}(\Omega,X)$.
\end{pf}

\begin{pf} [{\bf Proof of Theorem \ref{tthm-1.9-atm}}]
	We first prove the inequality $(1)$.  By going the similar line of arguments as in the proof of Theorem \ref{tthm-1.8-atm}, from \eqref{ee-1.5-atm}, we obtain 
	\begin{equation*}
		\sum_{m=0}^{\infty} \sum_{|\alpha|=m} (\norm{U(a_{\alpha})}^p_{Y} + \norm{U(b_{\alpha})}^p_{Y})(Rr)^{p \alpha} \leq \norm{U}^p \norm{a_{0}}^p +  \sum_{m=1}^{\infty} R^{pm} \lambda^{p} \norm{\sum_{|\alpha|=m} a_{\alpha}\, z^{\alpha} +  \sum_{|\alpha|=m} b^{*}_{\alpha}\, \bar{z}^{\alpha}}^p_{\Omega,X}.
	\end{equation*} 
	Now, Lemma \ref{lem-3.1} gives
	\begin{align*}
		\norm{\sum_{|\alpha|=m} a_{\alpha} z^{\alpha}}
		& \leq \norm{\sum_{|\alpha|=m} \left(\frac{a_{\alpha} + b_{\alpha}}{2}\right)z^{\alpha}} + \norm{\sum_{|\alpha|=m} \left(\frac{a_{\alpha} - b_{\alpha}}{2}\right)z^{\alpha}} \nonumber \\
		& \leq 4 \, \norm{\norm{f}_{\Omega,X}\,I-\real(a_{0})}.
	\end{align*}
	The last inequality remains valid if the left hand quantity replaced by $\norm{\sum_{|\alpha|=m} b_{\alpha} z^{\alpha}}$. Thus, 
	\begin{equation} \label{ee-1.7-atm}
		\norm{\sum_{|\alpha|=m} a_{\alpha}\, z^{\alpha} +  \sum_{|\alpha|=m} b^{*}_{\alpha}\, \bar{z}^{\alpha}}^p \leq 8^p \, \norm{\norm{f}_{\Omega,X}\,I-\real(a_{0})}^p.
	\end{equation}
	Using \eqref{ee-1.7-atm} and the fact $f(0)=a_0=\mu I$, $0 \leq \mu \leq \norm{f}_{\Omega,X}$, a simple computation shows that 
	\begin{equation*}
		\sum_{m=0}^{\infty} \sum_{|\alpha|=m} (\norm{U(a_{\alpha})}^p_{Y} + \norm{U(b_{\alpha})}^p_{Y})(Rr)^{p \alpha} \leq \lambda^p \mu^p + \lambda^p 8^p \left(\norm{f}^p_{\Omega,X} - \mu^p \right)\frac{R^p}{1-R^p} \leq \lambda^p \norm{f}^p_{\Omega,X},
	\end{equation*}
	provided 
	$$R \leq \frac{1}{(1+8^p)^{1/p}}.$$ This yields
	$$ \frac{1}{(1+8^p)^{1/p}}\,\frac{1}{n} \, \sum_{i=1}^{n} r_{i} \leq AP^{*}_{\lambda}(B_{Z}, p,U),$$
	and hence, the desired bound follows. 
	\vspace{1mm}
	
	By going the similar line of arguments as in the proof of part $(2)$ of Theorem \ref{tthm-1.8-atm}, using  \eqref{ee-1.6-atm} and \eqref{ee-1.7-atm}, we obtain
	\begin{equation*}
		\sum_{m=0}^{\infty} \sum_{|\alpha|=m} (\norm{U(a_{\alpha})}^p_{Y} + \norm{U(b_{\alpha})}^p_{Y})(Rr)^{p \alpha} \leq \lambda^p \mu^p + 8^p \left(\norm{f}^p_{\Omega,X} - \mu^p \right)\frac{R^p}{1-R^p} \leq \lambda^p \norm{f}^p_{\Omega,X},
	\end{equation*}
	provided 
	$$R \leq \frac{\lambda}{(\lambda^p+8^p)^{1/p}}.$$ This yields
	$$ \frac{\lambda}{(\lambda^p+8^p)^{1/p}}\,\frac{1}{n} \, \sum_{i=1}^{n} r_{i} \leq AP^{*}_{\lambda}(B_{Z}, p,U),$$
	and hence, the desired bound follows. The proof is now complete.
\end{pf}

\begin{pf}[{\bf Proof of Theorem \ref{thm-1.2-atm}}]
Let $X=\mathcal{B}(\mathcal{H})$. We shall first prove the lower bound of $AP_{\lambda}(B_{Z}, p,X)$. The following lower bound of $R_{\lambda}(B_{Z}, p,X)$ was established in \cite[Theorem 1.2]{Himadri-local-Banach-1} in the case where $X$ is finite dimensional:
\begin{equation*}
	R_{\lambda}(B_{Z}, p,X) \geq	\begin{cases}
		E_{1}(X)\,  \max \left\{\frac{1}{\sqrt{n}\, \norm{Id: \ell^n_{2} \rightarrow Z}}, \, \frac{1}{e\, \norm{Id: Z \rightarrow \ell^n_{1}}}\right\}\, \left(\frac{\lambda^p -1}{2\lambda^p - 1}\right)^{\frac{1}{p}} & \text{for $p=1$}, \\[2mm]
		E_{3}\norm{I:\ell^n _2 \rightarrow Z}^{-\frac{2}{p}}\, \left(\frac{\lambda^p -1}{2\lambda^p - 1}\right)^{\frac{1}{p}}\, & \text{for $p \geq 2$}, \\[2mm]
		
		E_{2}(X)\,  \max \left\{\frac{1}{(\sqrt{n})^{1-\theta} \norm{Id: \ell^n_{2} \rightarrow Z}},  \frac{\norm{I:\ell^n _2 \rightarrow Z}^{-\theta}}{e^{(1-\theta)} \norm{Id: Z \rightarrow \ell^n_{1}}^{1-\theta}}\right\}\left(\frac{\lambda^p -1}{2\lambda^p - 1}\right)^{\frac{1}{p}} \hspace{-2mm} & \text{for $1<p< 2$}, 
	\end{cases}
\end{equation*}
where $\theta = (2(p-1))/p$. Here, $E_{1}(X)$, $E_{2}(X)$, and $E_{3}$ are positive constants: $E_{1}(X)$ depends only on $X$, $E_{2}(X)$ depends on both $X$ and $p$, while $E_{3}$ is independent of $X$ and depends only on $p$. Thus, by combining Theorem \ref{thm-1.8-atm} with the above estimate of $R_{\lambda}(B_{Z}, p,X)$, the desired lower bound of $AP_{\lambda}(B_{Z}, p,X)$ follows. 
\vspace{1mm}

We now proceed to estimate the upper bound. In this regard, we first obtain a lower bound of $AP_{\lambda}(B_{Z}, p,X)$ when $Z=\ell^n_q$, $1 \leq q \leq \infty$, and then the desired lower bound of $AP_{\lambda}(B_{Z}, p,X)$ for general Banach space $Z$ follows from Lemma \ref{lem-3.6}. Let $r=(r_1,.\,.\,., r_n)\in \mathbb{R}^n_{\geq 0}$ be such that 
\begin{equation*}
	\sum_{|\alpha|=m} (\norm{c_{\alpha}}^p + \norm{d_{\alpha}}^p)r^{p \alpha} \leq \lambda^p \, \norm{g}^p_{B_{\ell^n_q}, X}	
\end{equation*}
for all $g(z)=\sum_{|\alpha|=m} c_{\alpha}\, z^{\alpha} +  \sum_{|\alpha|=m} d^{*}_{\alpha}\, \bar{z}^{\alpha}\in \mathcal{PH}(^mB_{\ell^n_q},X)$. Then, for $\varepsilon_{\alpha} \in \{-1,1\}$ and $I$ the identity operator on $\mathcal{H}$, the $m$-homogeneous polynomial
\begin{equation*}
f(z)=\sum_{|\alpha|=m} \varepsilon_{\alpha}\, \frac{m!}{\alpha!}z^{\alpha}\, I +  \sum_{|\alpha|=m} \varepsilon_{\alpha}\,\frac{m!}{\alpha!}\, \bar{z}^{\alpha} I, \,\, z \in B_{\ell^n_q},
\end{equation*} 
satisfies
\begin{equation} \label{ee-5.4-added}
2\, \sum_{|\alpha|=m}  \left(\frac{m!}{\alpha!}\right)^p\, r^{p\alpha} \leq \lambda^p \,\norm{f}^p_{B_{\ell^n_q},X}.
\end{equation}
Using the facts $\norm{r}_1 \leq n^{1-\frac{1}{p}}\, \norm{r}_p$ and $m! /\alpha! \geq 1$, a simple computation shows that
\begin{align*}
\left(\sum_{i=1}^{n}r_i\right)^{pm} 
& \leq \left(n^{1-\frac{1}{p}}\right)^{pm}\, \left(\sum_{i=1}^{n}r^p_i\right)^m \nonumber \\
& = \left(n^{1-\frac{1}{p}}\right)^{pm}\, \sum_{|\alpha|=m} \frac{m!}{\alpha!}\, r^{p \alpha} \nonumber \\
& \leq \left(n^{1-\frac{1}{p}}\right)^{pm}\, \sum_{|\alpha|=m} \left(\frac{m!}{\alpha!}\right)^p\, r^{p \alpha}.
\end{align*}
In view of the preceding inequality, from \eqref{ee-5.4-added}, we deduce that
\begin{align} \label{ee-5.5-added}
2\, \left(\sum_{i=1}^{n}r_i\right)^{pm} 
& \leq \lambda^p\, \left(n^{1-\frac{1}{p}}\right)^{pm}\, \norm{\sum_{|\alpha|=m} \varepsilon_{\alpha}\, \frac{m!}{\alpha!}z^{\alpha}\, I +  \sum_{|\alpha|=m} \varepsilon_{\alpha}\,\frac{m!}{\alpha!}\, \bar{z}^{\alpha} I}^p_{B_{\ell^n_q},X} \nonumber\\
& \leq \lambda^p\, \left(n^{1-\frac{1}{p}}\right)^{pm}\, \left(\norm{\sum_{|\alpha|=m} \varepsilon_{\alpha}\, \frac{m!}{\alpha!}z^{\alpha}\, I}_{B_{\ell^n_q},X} + \norm{\sum_{|\alpha|=m} \varepsilon_{\alpha}\,\frac{m!}{\alpha!}\, \bar{z}^{\alpha} I}_{B_{\ell^n_q},X}\right)^p \nonumber \\
& = \lambda^p\, \left(n^{1-\frac{1}{p}}\right)^{pm}\, \left(\norm{\sum_{|\alpha|=m} \varepsilon_{\alpha}\, \frac{m!}{\alpha!}z^{\alpha}}_{B_{\ell^n_q},\mathbb{C}} + \norm{\sum_{|\alpha|=m} \varepsilon_{\alpha}\,\frac{m!}{\alpha!}\, \bar{z}^{\alpha}}_{B_{\ell^n_q},\mathbb{C}}\right)^p \nonumber\\
&= \lambda^p\, \left(n^{1-\frac{1}{p}}\right)^{pm}\, \left(\norm{\sum_{|\alpha|=m} \varepsilon_{\alpha}\, \frac{m!}{\alpha!}z^{\alpha}}_{B_{\ell^n_q},\mathbb{C}} + \norm{\sum_{|\alpha|=m} \varepsilon_{\alpha}\,\frac{m!}{\alpha!}\, z^{\alpha} }_{B_{\ell^n_q},\mathbb{C}}\right)^p \nonumber\\
&= 2^p\, \lambda^p\, \left(n^{1-\frac{1}{p}}\right)^{pm}\,\norm{\sum_{|\alpha|=m} \varepsilon_{\alpha}\,\frac{m!}{\alpha!}\, z^{\alpha} }^p_{B_{\ell^n_q},\mathbb{C}}.
\end{align}
We now recall the following probabilistic estimate due to Boas \cite[Theorem 4]{boas-2000}:
\begin{equation*}
\norm{\sum_{|\alpha|=m} \varepsilon_{\alpha}\,\frac{m!}{\alpha!}\, z^{\alpha} }_{B_{\ell^n_q},\mathbb{C}} \leq 
\gamma(m,n,q):= \sqrt{32m\, \log(6m)}
\begin{cases}
n^{\frac{1}{2}}\, (m!)^{1-\frac{1}{q}}\, & \mbox{if}\,\, q \leq 2 \\
n^{\frac{1}{2}+ \left(\frac{1}{2}-\frac{1}{q}\right)m}\, (m!)^{\frac{1}{2}}\, & \mbox{if}\,\, q \geq 2.
\end{cases}
\end{equation*}
This estimate, together with \eqref{ee-5.5-added}, yields
\begin{equation*}
	\frac{\sum_{i=1}^{n}r_i}{n} \leq \left(2^{1-\frac{1}{p}}\right)^m \lambda^{\frac{1}{m}}\, n^{-\frac{1}{p}}\, (\gamma(m,n,q))^{\frac{1}{m}}.
\end{equation*}
Consequently,
\begin{align} \label{ee-5.6-added}
AP^m_{\lambda}(B_{\ell^n_q}, p,X) 
& \leq \left(2^{1-\frac{1}{p}}\right)^m \lambda^{\frac{1}{m}}\, n^{-\frac{1}{p}}\, (\gamma(m,n,q))^{\frac{1}{m}} \nonumber\\
& \leq C_{1}\,\lambda^{\frac{1}{m}}\, \left(m\, \log\,m (m!)^{2 \left(1-\frac{1}{\min\{2,q\}}\right)}\,n\right)^{\frac{1}{2m}}\, \left(\frac{1}{n}\right)^{\frac{1}{p}+\frac{1}{\max\{2,q\}}-\frac{1}{2}} 
\end{align}
for $C_1>0$, a uniform constant. We already know that $AP_{\lambda}(B_{\ell^n_q}, p,X) \leq AP^m_{\lambda}(B_{\ell^n_q}, p,X)$, so for the desired upper bound of $AP_{\lambda}(B_{\ell^n_q}, p,X) $, it is enough to estimate the right side quantity of  \eqref{ee-5.6-added}. Write $t=1-1/\min\{2,q\}$ for short. By Stirling's formula, we have $m! \sim \sqrt{2\, \pi m} (m/e)^m$. Hence
\begin{equation*}
(m!)^{2t} = \left(\frac{m}{e}\right)^{2tm} (2 \pi m)^{t} \left(1+O\left(\frac{1}{m}\right)\right),
\end{equation*}
and so the right side quantity of \eqref{ee-5.6-added} becomes less than or equals to
\begin{equation} \label{ee-5.7-added}
C_1 \lambda^{\frac{1}{m}}\, (m\, \log\,m \, n)^{\frac{1}{2m}}\, \left(\frac{m}{e}\right)^{t} (2 \pi m)^{\frac{t}{2m}} \left(1+O\left(\frac{1}{m}\right)\right)^{\frac{1}{2m}}\, \, \left(\frac{1}{n}\right)^{\frac{1}{p}+\frac{1}{\max\{2,q\}}-\frac{1}{2}}.
\end{equation}
Since $(2 \pi m)^{\frac{t}{2m}} \rightarrow 1$ as $m \rightarrow \infty$, the sequence $\{(2 \pi m)^{\frac{t}{2m}}\}$ is bounded. On the other hand, $(1+O(1/m))^{\frac{1}{2m}}$ is also bounded by uniform constant. We now want to estimate the quantity $(m\, \log\,m \, n)^{\frac{1}{2m}}$. Take $m=[\log\,n]$. Then
\begin{equation*}
\log\, (m\, \log\,m \, n)=\log\, n + \log\, m +\log\log\,m= \log \,n + O(\log\log \,n).
\end{equation*}
Hence,
\begin{equation*}
(m\, \log\,m \, n)^{\frac{1}{2m}}= \exp \left(\frac{\log \,n + O(\log\log \,n)}{2m}\right).
\end{equation*}
Since $m \sim \log \, n$,
\begin{equation*}
\frac{\log \,n}{2m} \leq \frac{1}{2} + o(1)\,\, \mbox{and} \, \, O\left(\frac{\log\log \,n}{\log \, n}\right) \rightarrow 0.
\end{equation*}
Therefore,
\begin{equation*}
(m\, \log\,m \, n)^{\frac{1}{2m}} \leq C_2	
\end{equation*}
for some constant $C_2>0$. Combining these estimates, the quantity in \eqref{ee-5.7-added} becomes less than or equals to
\begin{equation*}
C_{3}\, \lambda^{\frac{1}{m}}\, m^t\,	\left(\frac{1}{n}\right)^{\frac{1}{p}+\frac{1}{\max\{2,q\}}-\frac{1}{2}}, 
\end{equation*}
which finally less than or equals to
\begin{equation*}
C\, \lambda^{\frac{1}{\log\,n}}\,\, (\log\,n)^{1-\frac{1}{\min\{2,q\}}}\,\,\left(\frac{1}{n}\right)^{\frac{1}{p}+\frac{1}{\max\{2,q\}}-\frac{1}{2}}
\end{equation*}
for some universal constant $C>0$.
Thus,
\begin{equation} \label{ee-5.8-added}
AP_{\lambda}(B_{\ell^n_q}, p,X) \leq C\, \lambda^{\frac{1}{\log\,n}}\,\, (\log\,n)^{1-\frac{1}{\min\{2,q\}}}\,\,\left(\frac{1}{n}\right)^{\frac{1}{p}+\frac{1}{\max\{2,q\}}-\frac{1}{2}}.
\end{equation}
For the upper bound of $AP_{\lambda}(B_{Z}, p,X)$, by making use of Lemma \ref{lem-3.6} and \eqref{ee-5.8-added}, we establish
\begin{equation*}
AP_{\lambda}(B_{Z}, p,X) \leq C\, \lambda^{\frac{1}{\log\,n}}\, \norm{Id:Z \rightarrow \ell^n_q}\, (\log\,n)^{1-\frac{1}{\min\{2,q\}}}\,\,\left(\frac{1}{n}\right)^{\frac{1}{p}+\frac{1}{\max\{2,q\}}-\frac{1}{2}}.
\end{equation*}
This completes the proof.
\end{pf}
\begin{pf} [{\bf Proof of Corollary \ref{cor-1.4}}]
Let $X=\mathcal{B}(\mathcal{H})$. The following lower bound of $R_{\lambda}(B_{\ell^n_q}, p,X)$ was established in \cite[Corollary 1.2]{Himadri-local-Banach-1} in the case where $X$ is finite dimensional:
\begin{equation*}
	R_{\lambda}(B_{\ell^n _q}, p,X) \geq	\begin{cases}
		E'_{3}(X)\,\left(\frac{\lambda^p -1}{2\lambda^p - 1}\right)^{\frac{1}{p}} \, \left(\frac{\log n}{n}\right)^{ \left(1- \frac{1}{\min\{q,2\}}\right)} & \text{for $p=1$}, \\[2mm]
		E'_{2} \left(\frac{\lambda^p -1}{2\lambda^p - 1}\right)^{\frac{1}{p}}\, n^{-\frac{1}{p}} & \text{for $p \geq q$}, \\[2mm]
		
		E_{4}(X)\,  \left(\frac{\lambda^p -1}{2\lambda^p - 1}\right)^{\frac{1}{p}} \,n^{-\frac{p-1}{p(q-1)}}\, \left( \frac{\log n}{n}\right)^{\left(1- \frac{1}{\min\{q,2\}}\right)\, \frac{q-p}{p(q-1)}} & \text{for $1<p< q$}. 
	\end{cases}
\end{equation*}
Here, $E'_{3}(X)$, $E'_{2}$, and $E_{4}(X)$ are positive constants: $E'_{3}(X)$ depends only on $X$, $E_{4}(X)$ depends on both $X$ and $p$, while $E'_{2}$ is independent of $X$ and depends only on $p$. Thus, by combining Theorem \ref{thm-1.8-atm} with the preceding estimate of $R_{\lambda}(B_{\ell^n_q}, p,X)$, the desired lower bound of $AP_{\lambda}(B_{\ell^n_q}, p,X)$ follows. The proof is now complete.
\end{pf}
\begin{pf} [{\bf Proof of Theorem \ref{thm-1.3-a}}]
	We start by recalling two classical facts from the theory of finite-dimensional Banach spaces with unconditional bases. Let $W=(\mathbb{C}^n, ||.||)$ be a Banach space whose canonical unit vectors $e_{k}$ form a normalized $1$-unconditional basis. Then the operator norm of the identity embedding from $W$ into $\ell^n_2$ is given by $\norm{Id: W \rightarrow \ell_2^n}=\sup_{z \in B_{W}}(\sum_{m=1}^{n}|z_{m}|^2)^{1/2}$. On the other hand, the identity embedding from $Y$ into $\ell^n_1$ satisfies  $\norm{Id: W \rightarrow \ell_1^n}=\sup_{z \in B_{W}} \sum_{m=1}^{n}|z_{m}|=\norm{\sum_{m=1}^{n}e^{*}_{m}}_{W^{*}}$ for all $z \in W$. We are now in a position to establish the assertions of the theorem.
	\vspace{1mm}
	
	Since $Z \subset \ell_2$, we may assume that the inclusion map $Z \subset \ell_2$ has norm $1$. Under this normalization, the claim $(1)$ follows directly from Theorem \ref{thm-1.2-atm} together with the above identities for the identity maps.
	\vspace{1mm}
	
	We shall now prove part $(2)$. By a result of \cite[Theorem 1.d.5]{Lindenstrauss-book}, we may assume without loss of generality that the constant $M^{(2)}(X)=1$. This implies that $\norm{Id: \ell_2^n \rightarrow Z_{n}}=1$. We now recall the following important fact from  \cite[(5.3), p. 191]{defant-2003}. If $W=(\mathbb{C}^n, ||.||)$ be a symmetric Banach space such that $M^{(2)}(W)=1$, then 
	\begin{equation} \label{e-5.4-added}
	\norm{Id: W \rightarrow \ell^n_2}= \frac{\norm{\sum_{m=1}^{n}e^{*}_{m}}_{W^{*}}}{\sqrt{n}}.
	\end{equation}
	 The conclusion then follows by combining Theorem \ref{thm-1.2-atm}, the above facts for identity embeddings, and the inequality \eqref{e-5.4-added}. This concludes the proof.	
\end{pf}
\begin{pf} [{\bf Proof of Theorem \ref{thm-1.3-atm}}]
	Let $X=\mathcal{B}(\mathcal{H})$. The lower bound of $AP_{\lambda}(B_{\ell^n _q}, p,X)$ follows from Theorem \ref{thm-1.8-atm} by observing the fact 
	$$\norm{Id: B_{\ell^n_q} \rightarrow B_{\ell^n_1}}= n^{1-\frac{1}{q}}$$
	for any $1 \leq q \leq \infty$.
	Now, we want to find the upper bound of $AP_{\lambda}(B_{\ell^n _q}, p,X)$. Let $r=(r_1,.\,.\,., r_n)\in \mathbb{R}^n_{\geq 0}$ be such that 
	\begin{equation} \label{e-1.5-atm}
		\sum_{k=1}^{n} (\norm{c_{k}}^p + \norm{d_{k}}^p)r^{p} \leq  \norm{g}^p_{B_{Z}, \mathcal{B}(\mathcal{H})}	
	\end{equation}
	for all $g(z)=\sum_{k=1}^{n}  c_{k}\, z_{k} + \sum_{k=1}^{n} d^{*}_{k}\, \bar{z}_{k}\in \mathcal{PH}(^1B_{z},\mathcal{B}(\mathcal{H}))$. We shall now make use of Theorem \ref{thm-A}. By considering $z=e_{j}$ in Theorem \ref{thm-A}, for each $n \in \mathbb{N}$ there exist $a_{1}, \ldots, a_{n} \in X$ such that $1/2 \leq \norm{a_{j}} $ for each $1\leq j \leq n$. This inequality together with  H\"{o}lder's inequality, by considering the $1$-homogeneous polynomial $F \in \mathcal{PH}(^1B_{Z},X)$ of the form $F(z)=\sum_{j=1}^{n} z_{j}a_{j}+ \sum_{j=1}^{n} a^{*}_{j}\bar{z_{j}}$, yields
	\begin{align*} \label{e-1.39-atm}
		\sum_{k=1}^{n} r_{k}  \leq \sum_{k=1}^{n}(\norm{a_{j}}+\norm{a_{j}}) r_{k} \leq \left(\sum_{k=1}^{n} (\norm{a_{k}}+\norm{b_{k}})^p r^p_{k}\right)^{\frac{1}{p}} (\sum_{k=1}^{n} 1)^{1-\frac{1}{p}}  = 2^{1-\frac{1}{p}} n^{1 -\frac{1}{p}} \norm{F}_{B_{\ell^n_q},X}.  
	\end{align*} 
	Using Theorem \ref{thm-A}, from the last inequality we obtain 
	\begin{equation} \label{e-1.6-atm}
		\frac{1}{n} \, \sum_{k=1}^{n} r_{k} \leq 2^{1-\frac{1}{p}} n^{ -\frac{1}{p}} \sup_{z \in B_{\ell^n _q}} \norm{z}_{Cot(X)}.
	\end{equation}
	For $q \leq Cot(X)$, we have $\sup_{z \in B_{\ell^n _q}} \norm{z}_{Cot(X)} \leq 1$, which together with \eqref{e-1.6-atm} gives 
	$$AP^1(B_{\ell^n _q}, p,X)\leq 2^{1-\frac{1}{p}} n^{ -\frac{1}{p}}.$$ On the other hand, for $q>Cot(X)$, $$\sup_{z \in B_{\ell^n _q}} \norm{z}_{Cot(X)} \leq n^{1/Cot(X)\, - \,1/q},$$ and so from \eqref{e-1.6-atm} we obtain 
	$$AP^1(B_{\ell^n _q}, p,X) \leq 2^{1-\frac{1}{p}}\, n^{1/Cot(X)\, - \,1/q -1/p}.$$ It is known that $AP_{\lambda}(B_{\ell^n _q}, p,X)\leq  AP^1_{\lambda}(B_{\ell^n _q}, p,X)=\lambda\,AP^1(B_{\ell^n _q}, p,X)$, and hence the upper bound of $AP_{\lambda}(B_{\ell^n _q}, p,X)$ immediately follows. This completes the proof.
\end{pf} 
		
	\begin{pf} [{\bf Proof of Theorem \ref{thm-1.8}}]
	We shall make use of the following result due to \cite[Theorem 5.3]{defant-2012}. If $Y$ is a $s$-concave Banach lattice, with $2\leq s < \infty$, and $U: X \rightarrow Y$ an $(t,1)$-summing operator with $1 \leq t \leq s$. Then there is a constant $C>0$ such that
	\begin{equation} \label{e-5.7-added}
		\left(\sum_{|\alpha|=m} \norm{U(c_{\alpha})}^{\beta}_{Y}\right)^{\frac{1}{\beta}} \leq C^m \, \norm{Q}_{\mathbb{D}^n,X}
	\end{equation}
	holds for every homogeneous polynomial $Q: \mathbb{D}^n \rightarrow X$ with $Q(z)=\sum_{|\alpha|=m} c_{\alpha}\, z^{\alpha}$ on $\mathbb{D}^n$, where 
	\begin{equation*}
		\beta= \frac{mst}{s+ (m-1)t}.
	\end{equation*}	
The proof of the theorem relies on the definition of complete Reinhardt domains and identifying $B_z$ as such a domain.  Recall that a domain $\Omega \subset \mathbb{C}^n$ is called complete Reinhardt if $(z_{1}, \ldots,z_{n}) \in \Omega$ then the new tuple $(\xi_1 z_1, \ldots , \xi _nz_n) \in \Omega$ for each $\xi_j \in \overline{\mathbb{D}}$, $j=1, \ldots,n$.
Here, $Z=(\mathbb{C}^n, ||.||)$ is an $n$-dimensional Banach space for which the canonical basis vectors $e_{k}$ form a normalized $1$-unconditional basis. That is, equivalent to say the open unit ball $B_{Z}$ is a complete Reinhardt domain.
 
Fix $z \in B_{Z}$. For every $m$-homogeneous polynomial $P: B_{Z} \rightarrow X$ with $P(z)=\sum_{|\alpha|=m} a_{\alpha}\, z^{\alpha}$ on $B_{Z}$, we consider the following $m$-homogeneous polynomial $\tilde{Q}: \mathbb{D}^n \rightarrow X$ by 
\begin{equation*}
	\tilde{Q}(w_1, \ldots,w_n):=P(z_1w_1, \ldots , z_nw_n)= \sum_{|\alpha|=m} (a_{\alpha}\, z^{\alpha}) w^{\alpha}
\end{equation*}
for $w=(w_1, \ldots ,w_n) \in \mathbb{D}^n$. The map $\tilde{Q}$ is well defined, since $B_Z$ is a complete Reinhardt domain and hence the point $(z_1w_1, \ldots , z_nw_n)$ belongs to $B_Z$. Clearly, $\norm{\tilde{Q}}_{\mathbb{D}^n,X} \leq \norm{P}_{B_{Z},X}$. Now, by virtue of \eqref{e-5.7-added}, we obtain
\begin{equation*} 
	\left(\sum_{|\alpha|=m} \norm{U(a_{\alpha})z^{\alpha}}^{\beta}_{Y}\right)^{\frac{1}{\beta}} \leq C^m \, \norm{\tilde{Q}}_{\mathbb{D}^n,X} \leq C^m \, \norm{P}_{B_{Z},X}.
\end{equation*}
This completes the proof.
	\end{pf}
	
\begin{pf} [{\bf Proof of Theorem \ref{thm-1.9}}]
We first want to prove the part $(1)$. Let $X=\mathcal{B}(\mathcal{H}_1)$ and $Y=\mathcal{B}(\mathcal{H}_2)$. By the definition of Bohr radius, it not difficult to see that
\begin{equation} \label{ee-5.8}
R_{\lambda}(B_{Z}, p,U) \geq \max\left\{R_{\frac{\lambda}{\norm{U}}}(B_{Z}, p,X), \, \, R_{\frac{\lambda}{\norm{U}}}(B_{Z}, p,Y)\right\}.
\end{equation}
In view of \cite[Theorem 1.5]{Himadri-local-Banach-1} and \eqref{ee-5.8}, for $1\leq p<\infty$ and $1 \leq q \leq \infty $, we deduce that
\begin{align} \label{e-5.9-added}
R_{\lambda}(B_{Z}, p,U) \geq \frac{\Phi_{1}(n,\lambda,p,q)}{S(B_{Z},B_{\ell^n _q})S(B_{\ell^n _q},B_{Z})},
\end{align} 
where
\begin{equation*}
\Phi_{1}(n,\lambda,p,q):=	\begin{cases}
	
	E_{5}\,\dfrac{\left(\lambda^p -\norm{U}^p\right)^{\frac{1}{p}}}{\lambda} & \text{for $p>\min \{Cot(X),Cot(Y),q\}$}, \\[2mm]
	E_{6}\, \dfrac{\left(\lambda^p -\norm{U}^p\right)^{\frac{1}{p}}}{\lambda}\,\, n^{{\frac{1}{\min\{t,q\}}}-\frac{1}{p}} & \text{for $p\leq \min \{Cot(X),Cot(Y),q\}$}.
\end{cases}
\end{equation*}
 Here $E_{5}>0$ and $E_{6}>0$ are constants independent of $X$ and $n$. The lower bound of $AP_{\lambda}(B_{Z}, p,U)$ follows from Theorem \ref{thm-1.8-atm} and \eqref{e-5.9-added}.
 
 We shall now prove the part $(2)$. Now let $X=\mathcal{B}(\mathcal{H})$.
Let 
$$P(z)=\sum_{|\alpha|=m} a_{\alpha}\, z^{\alpha} + \sum_{|\beta|=m} b^{*}_{\alpha}\, \bar{z}^{\alpha} \in \mathcal{PH}(^m Z,X).$$ Consider the holomorphic polynomial $$H(z):=\sum_{|\alpha|=m} \left(\frac{a_{\alpha} + b_{\alpha}}{2}\right)z^{\alpha}$$ on $B_{Z}$. In view of H\"{o}lder's inequality and Theorem \ref{thm-1.8}, we obtain
\begin{align*}
&\sum_{|\alpha|=m} \norm{U\left(\frac{a_{\alpha} + b_{\alpha}}{2}\right)z^{\alpha}}^p \nonumber\\
&\leq  \left(\sum_{|\alpha|=m} \norm{U\left(\frac{a_{\alpha} + b_{\alpha}}{2}\right)z^{\alpha}}\right)^p \nonumber \\
& \leq \left[\bigg(\sum_{|\alpha|=m}1\bigg)^{\frac{(s-1)mt-s+t}{mst}} \times \left(\sum_{|\alpha|=m} \norm{U\left(\frac{a_{\alpha} + b_{\alpha}}{2}\right)z^{\alpha}}^{\frac{mst}{s+(m-1)t}}\right)^{\frac{s+(m-1)t}{mst}}\right]^p \nonumber \\
& \leq \bigg(\sum_{|\alpha|=m}1\bigg)^{\frac{((s-1)mt-s+t)p}{mst}} c_{1}^{mp}\, \norm{H}^p_{\Omega,X} \leq \bigg(\sum_{|\alpha|=m}1\bigg)^{\frac{((s-1)mt-s+t)p}{mst}} c_{1}^{mp}\, 2^p \, \norm{P}^p_{B_{Z},X}.
\end{align*}
Similarly, by considering the holomorphic polynomial $G(z):=\sum_{|\alpha|=m} \left(\frac{a_{\alpha} - b_{\alpha}}{2}\right)z^{\alpha}$ on $B_{Z}$, we deduce that
\begin{equation*}
\sum_{|\alpha|=m} \norm{U\left(\frac{a_{\alpha} - b_{\alpha}}{2}\right)z^{\alpha}}^p \leq \bigg(\sum_{|\alpha|=m}1\bigg)^{\frac{((s-1)mt-s+t)p}{mst}} c_{2}^{mp}\, 2^p \, \norm{P}^p_{B_{Z},X}.
\end{equation*} 
A simple computation from last two inequalities shows that
\begin{align} \label{e-5.8}
&\left(\sum_{|\alpha|=m} \norm{U\left(\frac{a_{\alpha} + b_{\alpha}}{2}\right)z^{\alpha}}^p\right)^{1/p}  \nonumber \\
& \leq \left(\sum_{|\alpha|=m} \norm{U\left(\frac{a_{\alpha} + b_{\alpha}}{2}\right)z^{\alpha}}^p\right)^{1/p} + \left(\sum_{|\alpha|=m} \norm{U\left(\frac{a_{\alpha} - b_{\alpha}}{2}\right)z^{\alpha}}^p\right)^{1/p} \nonumber \\
& \leq 2(c^m_{1}+c^m_2) \bigg(\sum_{|\alpha|=m}1\bigg)^{\frac{(s-1)mt-s+t}{mst}} \norm{P}_{B_{Z},X}
\end{align} 
and 
\begin{equation} \label{e-5.9}
\left(\sum_{|\alpha|=m} \norm{U\left(\frac{a_{\alpha} - b_{\alpha}}{2}\right)z^{\alpha}}^p\right)^{1/p} \leq 2(c^m_{1}+c^m_2) \bigg(\sum_{|\alpha|=m}1\bigg)^{\frac{(s-1)mt-s+t}{mst}} \norm{P}_{B_{Z},X}.
\end{equation}
In view of \eqref{e-5.8} and \eqref{e-5.9}, we obtain 
\begin{equation} \label{e-5.10}
\sum_{|\alpha|=m} \left(\norm{U(a_{\alpha})}^p + \norm{U(b_{\alpha})}^p\right) |z^{\alpha}|^p \leq C^{pm} \bigg(\sum_{|\alpha|=m}1\bigg)^{\frac{((s-1)mt-s+t)p}{mst}} \norm{P}^p_{B_{Z},X}
\end{equation}
for some constant $C>0$. Using the estimate 
$$\sum_{|\alpha|=m} 1 =\binom{n+m-1}{m} \leq e^m \left(1+\frac{n}{m}\right)^m,
$$
from \eqref{e-5.10}, we deduce that 
\begin{equation} \label{e-5.11}
R^m(B_Z,p,U) \geq \tilde{D} \left(1+\frac{n}{m}\right)^{-\frac{(s-1)mt-s+t}{mst}}
\end{equation}
for some constant $\tilde{D}>0$. Since $R^m_{\lambda}(B_Z,p,U)=\lambda^{1/m} \, R^m(B_Z,p,U)$, \eqref{e-5.11} gives
\begin{equation} \label{e-5.12}
R^m_{\lambda}(B_Z,p,U) \geq\lambda^{\frac{1}{m}} \, \tilde{D} \left(1+\frac{n}{m}\right)^{-\frac{(s-1)mt-s+t}{mst}}.
\end{equation}
We now observe the elementary fact if $a,b>0$ and $0<p<1$, then $(a+b)^p \leq 2^p\, \max\{a^p, b^p\}$. Clearly, $$0<\frac{(s-1)mt-s+t}{mst}<1.$$
Thus, we have
\begin{equation*}
\left(1+\frac{n}{m}\right)^{\frac{(s-1)mt-s+t}{mst}} \leq 2^{\frac{(s-1)mt-s+t}{mst}} \max \left\{1, \left(\frac{n}{m}\right)^{\frac{(s-1)mt-s+t}{mst}}\right\} \leq 2\, \max \left\{1, \left(\frac{n}{m}\right)^{\frac{(s-1)m-(\frac{s}{t}-1)}{ms}}\right\}.
\end{equation*}
Consequently,
\begin{equation} \label{e-5.13}
	\left(1+\frac{n}{m}\right)^{-\frac{(s-1)mt-s+t}{mst}} \geq \frac{1}{2}\, \min \left\{1, \left(\frac{n}{m}\right)^{-\frac{(s-1)m-(\frac{s}{t}-1)}{ms}}\right\}=\frac{1}{2}\, \min \left\{1, \frac{m^{\frac{s-1}{s}}n^{\frac{(\frac{s}{t}-1)}{sm}}}{n^{\frac{s-1}{s}}m^{\frac{(\frac{s}{t}-1)}{sm}}}\right\}.
\end{equation}
It is not difficult  to see that the function $\psi: (0,\infty) \rightarrow \mathbb{R}$ defined by $$\psi(x)=x^{\frac{s-1}{s}} n^{\frac{(\frac{s}{t}-1)}{sm}},$$ attains a strict minimum at $x=\log \, n$ with the minimum value $$e^{\frac{(\frac{s}{t}-1)}{s}}\, (\log \, n)^{\frac{s-1}{s}}.$$ Finally, this fact together with Lemma \ref{lem-3.3} (1), \eqref{e-5.12}, and \eqref{e-5.13}, yields
\begin{equation} \label{e-5.15-added}
R_{\lambda}(B_Z,p,U) \geq D\,  \left(\frac{\lambda^p - \norm{U}^p}{2\lambda^p - \norm{U}^p}\right)^{\frac{1}{p}} \, \left(\frac{\log \, n}{n}\right)^{1-\frac{1}{s}}
\end{equation}
for some constant $D>0$. The lower bound of $AP_{\lambda}(B_{Z}, p,U)$ is a direct consequence Theorem \ref{thm-1.8-atm} and \eqref{e-5.15-added}.This completes the proof.
\comment{Observe that $1<t \leq s$, and so $(s-t)/mst \geq 0$. Using this fact, from \eqref{e-5.12}, we obtain 
	\begin{equation*}
		R^m_{\lambda}(B_Z,p,U) \geq\lambda^{\frac{1}{m}} \, \tilde{D} \left(1+\frac{n}{m}\right)^{-\frac{(s-1)}{s}}\, \left(1+\frac{n}{m}\right)^{\frac{(s-t)}{mst}} \geq \lambda^{\frac{1}{m}} \, \tilde{D}\, \left(1+\frac{n}{m}\right)^{-\frac{(s-1)}{s}}.
\end{equation*}}
\end{pf}

\begin{pf} [{\bf Proof of Theorem \ref{tthm-1.10-atm}}]
	Let $r=R_\lambda(\mathbb{D},p,U)$ and $f \in \mathcal{PH}(B_{z},X)$ be of the form \eqref{e-1.3-a}, where $X=\mathcal{B}(\mathcal{H})$.  Consider the function $\tilde{f}(z)=f(z, 0, \ldots,0)$ for $z \in \mathbb{D}$. Clearly, $\tilde{f}\in \mathcal{PH}(\mathbb{D},X)$ with
	the following expansion
	\begin{equation*}
		\tilde{f}(z)= \sum_{m=0}^{\infty} a_{(m,0,\ldots,0)} z^m + \sum_{m=1}^{\infty} b^{*}_{(m,0,\ldots,0)} \overline{z}^m, \,\,z \in \mathbb{D},
	\end{equation*}
	and $\norm{\tilde{f}}_{\mathbb{D},X}=\norm{f}_{B_{Z},X}$.
	Therefore, we have 
	\begin{align*}
		\sum_{m=0}^{\infty} \sum_{|\alpha|=m} (\norm{U(a_{\alpha})}^p_{Y} + \norm{U(b_{\alpha})}^p_{Y}) (r,0,\,.\,.\,.,0)^{p\alpha}
		&=\sum_{m=0}^{\infty}(\norm{U(a_{(m,0,\ldots,0)})}^p_{Y} + \norm{U(b_{(m,0,\ldots,0)})}^p_{Y}) r^{p\alpha}\nonumber,
	\end{align*}
	which is less than or equals to $\lambda^p \norm{\tilde{f}}^p_{\mathbb{D},X}=\lambda^p \norm{f}^p_{B_{Z},X}$.
	This shows that 
	$$\frac{r}{n}\leq AP_\lambda(B_{Z},p, U),$$
	 and hence we obtain $$\frac{R_\lambda(\mathbb{D},p,U)}{n} \leq AP_\lambda(B_{Z},p, U).
	 $$  Conversely, we prove that 
	\begin{equation*}
		AP_\lambda(B_{Z},p, U)\leq \norm{Id: Z \rightarrow \ell^n_{1}}\, \frac{R_\lambda(\mathbb{D},p,U)}{n^{1/p}}.
	\end{equation*}
	Suppose $r\in \mathbb{R}^n_{\geq 0}$ such that for all $f \in \mathcal{PH}(B_{z},X)$ be of the form \eqref{e-1.3-a},
	\begin{equation} \label{ee-1.8-atm}
		\sum_{m=0}^{\infty}	\sum_{|\alpha|=m} (\norm{U(a_{\alpha})}^p_{Y} + \norm{U(b_{\alpha})}^p_{Y})r^{p \alpha} \leq \lambda^{p} \norm{f}^p_{\Omega, X}.	
	\end{equation}
	We want to prove that 
	\begin{equation*}
		\frac{\sum_{j=1}^{n}r_j}{n}\leq \norm{Id: Z \rightarrow \ell^n_{1}}\, \frac{R_\lambda(\mathbb{D},p,U)}{n^{1/p}}.
	\end{equation*}
	Let $\tilde{f}:\mathbb{D}\rightarrow X \in \mathcal{PH}(\mathbb{D},X)$ be a bounded holomorphic map of the form $$\tilde{f}(\xi)=\sum_{m=0}^{\infty} c_{m} \xi^m + \sum_{m=1}^{\infty} d^{*}_m \overline{\xi}^m
	$$
	 for $\xi \in \mathbb{D}$. We now consider the function $\phi: B_{Z}\rightarrow \mathbb{C}$ defined by 
	\begin{equation} \label{e-function}
		\phi(z)=\frac{z_1+\cdots+z_n}{t},\quad z\in B_{Z},
	\end{equation}
	where $t= \norm{Id: Z \rightarrow \ell^n_{1}}$. Clearly, $\phi$ is holomorphic on $B_{Z}$ and $|\phi(z)|<1$ for all $z \in B_{Z}$.
	Now if we set $F=\tilde{f} \circ \phi $, then $F \in \mathcal{PH}(\Omega,X)$ with the following expansion 
	\begin{equation*}
		F(z)=\sum_{m=0}^{\infty} \sum_{|\alpha|=m} \left(\frac{c_{m}}{t^m} \frac{m!}{\alpha!}\right) z^{\alpha} +  \sum_{m=1}^{\infty} \sum_{|\alpha|=m} \left(\frac{d^{*}_{m}}{t^m} \frac{m!}{\alpha!}\right) \overline{z}^{\alpha},\,\, \,\, z \in B_{Z}
	\end{equation*}
	and $\norm{F}_{B_{Z},X}=\norm{\tilde{f}}_{\mathbb{D},X}$.
	By observing the fact $m!/\alpha! \geq 1$, for all $z\in B_{Z}$, we have
	\begin{align*}
		&\sum_{m=0}^{\infty} \sum_{|\alpha|=m} \left(\norm{U\left(\frac{c_{m}}{t^m} \frac{m!}{\alpha!}\right)}^p_{Y} + \norm{U\left(\frac{d_{m}}{t^m} \frac{m!}{\alpha!}\right)}^p_{Y}\right)|z|^{p \alpha} \nonumber \\
		& \geq \sum_{m=0}^{\infty} \sum_{|\alpha|=m} \left(\frac{m!}{\alpha!}\norm{U\left(\frac{c_{m}}{t^m} \right)}^p_{Y} + \frac{m!}{\alpha!}\norm{U\left(\frac{d_{m}}{t^m} \right)}^p_{Y}\right)|z|^{p \alpha} \nonumber \\
		& = \sum_{m=0}^{\infty} (\norm{U(c_{m})}^p_{Y} + \norm{U(d_{m})}^p_{Y})\, \left(\frac{\norm{z}_{p}}{t}\right)^{p m}.
	\end{align*}
	An observation, using preceding inequality and \eqref{ee-1.8-atm}, shows that
	\begin{align*}
		\sum_{m=0}^{\infty} (\norm{U(c_{m})}^p_{Y} + \norm{U(d_{m})}^p_{Y})\, \left(\frac{\norm{r}_{p}}{t}\right)^{p m}
		&\leq \sum_{m=0}^{\infty} \sum_{|\alpha|=m} \left(\norm{U\left(\frac{c_{m}}{t^m} \frac{m!}{\alpha!}\right)}^p_{Y} + \norm{U\left(\frac{d_{m}}{t^m} \frac{m!}{\alpha!}\right)}^p_{Y}\right)r^{p \alpha} \nonumber \\
		&\leq \lambda^p  \norm{F}^p_{B_{Z},X}=\lambda^p \norm{\tilde{f}}^p_{\mathbb{D},X},
	\end{align*}
	which implies that 
	$$
	\frac{\norm{r}_{p}}{t} \leq R_\lambda(\mathbb{D},p,U).
	$$ Again we have $\norm{r}_{1}\leq n^{(1-1/p)}\norm{r}_{p}$. Hence we obtain $$
	\frac{\norm{r}_{1}}{n}\leq \frac{t R_\lambda(\mathbb{D},p,U)}{n^{1/p}}.
	$$ 
	This completes the proof.
\end{pf}

\begin{pf} [{\bf Proof of Lemma \ref{lem-1.1}}]
Since $U$ is the identity operator on $\mathbb{C}$, we consider the complex valued harmonic functions on $\mathbb{D}$. Let $f:\mathbb{D} \rightarrow \mathbb{C}$ be a bounded harmonic function of the form
\begin{equation*}
	f(z)= \sum_{n=0}^{\infty} a_n\, z^n + \overline{\sum_{n=1}^{\infty} b_n\, z^n}, \,\,\, z \in \mathbb{D}.
\end{equation*}
Then, by \cite[Lemma 4]{Abu-2010}, we have 
\begin{equation*}
	|a_n|+ |b_n| \leq \frac{4}{\pi}\,\norm{f}_{\mathbb{D},\mathbb{C}}\,\,\,\,\, \mbox{for}\,\,\, n \geq 1.
\end{equation*}
By virtue of this estimate, for $|z|=r<1$, we obtain
\begin{align*}
|a_0|+ \sum_{n=1}^{\infty} (|a_n|+|b_n|)r^n 
& \leq |a_0|+ \frac{4}{\pi}\, \frac{r}{1-r}\,\norm{f}_{\mathbb{D},\mathbb{C}} \\
& \leq \norm{f}_{\mathbb{D},\mathbb{C}} + \frac{4}{\pi}\, \frac{r}{1-r}\,\norm{f}_{\mathbb{D},\mathbb{C}} \\
&= \left(1 + \frac{4}{\pi}\, \frac{r}{1-r} \right)\, \norm{f}_{\mathbb{D},\mathbb{C}} \leq \lambda \, \norm{f}_{\mathbb{D},\mathbb{C}},
\end{align*}
provided 
$$
r \leq \frac{(\lambda-1)\pi}{4+(\lambda-1)\pi}.
$$
Hence, the desired conclusion follows. The proof is now complete.
\comment{\\
	We now establish the part $(2)$. In this case, $U$ is the identity operator on $\mathcal{B}(\mathcal{H})$ and we consider the $\mathcal{B}(\mathcal{H})$ valued harmonic functions on $\mathbb{D}$. Let $f:\mathbb{D} \rightarrow \mathcal{B}(\mathcal{H})$ be a bounded harmonic function of the form
	\begin{equation*}
		f(z)= \sum_{n=0}^{\infty} a_n\, z^n + \sum_{n=1}^{\infty} b^{*}_n\, \overline{z}^n, \,\,\, z \in \mathbb{D}.
	\end{equation*}
	Then, by \cite[Lemma 4]{bhowmik-2021}, we have 
	\begin{equation*}
		||a_n||+ ||b_n|| \leq \frac{4}{\pi}\,\norm{f}_{\mathbb{D},\mathcal{B}(\mathcal{H})}\,\,\,\,\, \mbox{for}\,\,\, n \geq 1.
\end{equation*}}
\end{pf}

\noindent{\bf Statements and Declarations:}\\

\noindent{\bf Acknowledgment.} The author gratefully acknowledges the support received during the initial stages of this work while serving as a Senior Project Associate at IIT Bombay, through the \textit{Core Research Grant (CRG)}, sanctioned to Professor Sourav Pal by the Science and Engineering Research Board (SERB), Government of India. The author also sincerely acknowledges Professor B. V. Rajarama Bhat, ISI Bangalore, for support through the J. C. Bose Fellowship during the final stages of this manuscript. The author is currently supported by an NBHM Postdoctoral Fellowship from the Department of Atomic Energy (DAE), Government of India.\\

\noindent{\bf Conflict of interest.} The author declares that there is no conflict of interest regarding the publication of this paper.\\

\noindent{\bf Data availability statement.} Data sharing not applicable to this article as no datasets were generated or analysed during the current study.\\

\noindent{\bf Competing Interests.} The author declares none.\\

\noindent{\bf ORCID Id.} 0000-0001-9466-
178X

\end{document}